\documentclass{gtmon_a}
\pdfoutput=1

\usepackage{amscd}

\proceedingstitle{The interaction of finite-type and Gromov--Witten
invariants (BIRS 2003)}
\conferencestart{15 November 2003}
\conferenceend{20 November 2003}
\conferencename{The interaction of finite-type and Gromov--Witten
invariants}
\conferencelocation{Banff International Research Station, Banff, Alberta,
Canada}

\editor{David Auckly}
\givenname{David}
\surname{Auckly}

\editor{Jim Bryan}
\givenname{Jim}
\surname{Bryan}

\title[Enumerative geometry of stable maps]
{Enumerative geometry of stable maps with Lagrangian boundary conditions
and multiple covers of the disc}

\author{Sheldon Katz}
\givenname{Sheldon}
\surname{Katz}
\address{Departments of Mathematics and Physics\\
University of Illinois at Urbana-Champaign\\\newline
Urbana\\
Illinois 61801\\USA}
\address{Department of Mathematics\\
Oklahoma State University\\\newline
Stillwater\\
Oklahoma 74078\\USA}
\email{katz@math.uiuc.edu}

\author{Chiu-Chu Melissa Liu}
\givenname{Chiu-Chu Melissa}
\surname{Liu}
\address{Department of Mathematics\\Harvard University\\\newline
Cambridge\\
Massachusetts 02138\\USA}
\email{ccliu@math.harvard.edu}

\volumenumber{8}
\issuenumber{}
\publicationyear{2006}
\papernumber{1}
\startpage{1}
\endpage{47}

\doi{}
\MR{}
\Zbl{}


\received{22 January 2002}
\revised{}
\accepted{}
\published{22 April 2006 (in these proceedings)}
\publishedonline{22 April 2006 (in these proceedings)}
\proposed{}
\seconded{}
\corresponding{}
\version{}

\arxivreference{math.AG/0103074}

\dedicatory{Reproduced by kind permission of International Press
from:\newline 
{\rm Advances in Theoretical and Mathematical Physics, Volume 5 (2002) pages 1--49}}

\AtBeginDocument{\let\bar\wbar\let\tilde\wtilde\let\hat\what
\def\notin{\not\in}}
\allowdisplaybreaks

\makeatletter
\def\cnewtheorem#1[#2]#3{\newtheorem{#1}{#3}[subsection]
\expandafter\let\csname c@#1\endcsname\c@thm}
\def\dnewtheorem#1[#2]#3{\newtheorem{#1}{#3}[section]
\expandafter\let\csname c@#1\endcsname\c@tmm}

\newtheorem{thm}{Theorem}[subsection]
\cnewtheorem{pro}[thm]{Proposition}
\newtheorem{tmm}{Theorem}[section]
\dnewtheorem{propo}[tmm]{Proposition}
\cnewtheorem{cor}[thm]{Corollary}
\cnewtheorem{lm}[thm]{Lemma}
\theoremstyle{definition}
\cnewtheorem{example}[thm]{Example}
\cnewtheorem{rem}[thm]{Remark}
\cnewtheorem{df}[thm]{Definition}
\makeatother  

\makeautorefname{df}{Definition}

\newcommand{\overM}{{}\mskip3mu\overline{\mskip-3mu M\mskip-1mu}\mskip1mu}
\newcommand{\overSi}{{}\mskip2mu\overline{\mskip-1mu \Sigma\mskip-1mu}\mskip1mu}

\renewcommand{\Bar}{\overline}
\newcommand{\Cstar}{\mathbf{C}^*}
\newcommand{\N}{\mathcal{N}}
\newcommand{\CP}{\mathbf{P}}
\newcommand{\PP}{\mathbf{P}^1}
\newcommand{\MCY}{\overM_{g,0}(X,d\beta)}
\newcommand{\MPP}{\overM_{g,0}(\mathbf{P}^1, d)}
\newcommand{\MD}{\overM_{g;h}(D^2, S^1\mid d; n_1,\ldots, n_h)}
\newcommand{\MXL}{\overM_{g;h}(X,L\mid d\beta; n_1\gamma,\ldots, n_h\gamma)}
\newcommand{\MX}{\overM_{g;h}(X,L\mid \beta; \gamma_1,\ldots, \gamma_h)}
\newcommand{\fix}{F_{g;h|d;n_1,...,n_h}}
\newcommand{\OO}{\mathcal{O}_{\mathbf{P}^1}(-1)\oplus
                 \mathcal{O}_{\mathbf{P}^1}(-1)}
\newcommand{\OOd}{\mathcal{O}_{\mathbf{P}^1}(-d)\oplus
                 \mathcal{O}_{\mathbf{P}^1}(-d)}
\newcommand{\Oh}{\mathcal{O}}
\newcommand{\U}{\mathcal{U}}
\newcommand{\bS}{\partial\Sigma}
\newcommand{\inS}{\Sigma^\circ}
\newcommand{\Si}{\Sigma}
\newcommand{\si}{\sigma}
\newcommand{\tsi}{\tilde{\sigma}}
\newcommand{\hSi}{\hat{\Sigma}}
\newcommand{\hsi}{\hat{\sigma}}
\newcommand{\E}{\mathcal{E}}
\newcommand{\T}{\mathcal{T}}
\newcommand{\hodge}{\mathbb{E}}
\newcommand{\bra}{\langle}
\newcommand{\ket}{\rangle}

\begin{document}

\begin{asciiabstract}
In this paper, we present foundational material towards the development
of a rigorous enumerative theory of stable maps with Lagrangian
boundary conditions, ie stable maps from bordered Riemann surfaces
to a symplectic manifold, such that the boundary maps to a Lagrangian
submanifold. Our main application is to a situation where our proposed
theory leads to a well-defined algebro-geometric computation very
similar to well-known localization techniques in Gromov--Witten theory. In
particular, our computation of the invariants for multiple covers of a
generic disc bounding a special Lagrangian submanifold in a Calabi--Yau
threefold agrees completely with the original predictions of Ooguri
and Vafa based on string duality. Our proposed invariants depend more
generally on a discrete parameter which came to light in the work of
Aganagic, Klemm, and Vafa which was also based on duality, and our more
general calculations agree with theirs up to sign.
\end{asciiabstract}

\begin{abstract} 
In this paper, we present foundational material towards the development
of a rigorous enumerative theory of stable maps with Lagrangian
boundary conditions, ie stable maps from bordered Riemann surfaces
to a symplectic manifold, such that the boundary maps to a Lagrangian
submanifold. Our main application is to a situation where our proposed
theory leads to a well-defined algebro-geometric computation very
similar to well-known localization techniques in Gromov--Witten theory. In
particular, our computation of the invariants for multiple covers of a
generic disc bounding a special Lagrangian submanifold in a Calabi--Yau
threefold agrees completely with the original predictions of Ooguri
and Vafa based on string duality. Our proposed invariants depend more
generally on a discrete parameter which came to light in the work of
Aganagic, Klemm, and Vafa which was also based on duality, and our more
general calculations agree with theirs up to sign.
\end{abstract}

\maketitle

\section{Introduction}
\label{intro}
  
In recent years, we have witnessed the enormous impact that
theoretical physics has had on the geometry of Calabi--Yau manifolds.
One takes a version of string theory (or M--theory, or F--theory) and
compactifies the theory on a Calabi--Yau manifold, obtaining an
effective physical theory.  One of the key tools is {\em duality\/},
whose assertion is that several a priori distinct theories are related
in non-obvious ways.  Perhaps the best known duality is mirror
symmetry, which can be stated as the equality of the compactification
of type IIA string theory on a Calabi--Yau manifold $X$ with the
compactification of type IIB string theory on a mirror Calabi--Yau
manifold $X^\circ$.  This has deep mathematical consequences, which
are still far from being completely understood.  But there has
certainly been remarkable progress, in particular in the realm of
applications to enumerative geometry, as mirror symmetry (and
considerations of the topological string) inspired the rigorous
mathematical notions of the Gromov--Witten invariant and quantum
cohomology.  This has provided a precise context for the
string-theoretic prediction of the number of rational curves of
arbitrary degree on a quintic threefold (Candelas--de~la
Ossa--Green--Parkes~\cite{CdGP}).  Once the mathematical formulation
was clarified, then a mirror theorem could be formulated and proven
(Givental~\cite{G}, Lian--Liu--Yau~\cite{LLY}).  These aspects of the
development of mirror symmetry are discussed in detail in Cox and
Katz~\cite{CK}.

The enumerative consequences just discussed arose from considerations
of closed strings.  In closed string theory, closed strings sweep out
compact Riemann surfaces as the string moves in time.  Gromov--Witten
theory assigns invariants to spaces of maps from these Riemann
surfaces and their degenerations to collections of cohomology classes
on $X$.  

Recently, string theorists have obtained a better understanding of
the open string sector of string theory and how it relates to
dualities. Open strings sweep out compact Riemann surfaces with
boundary, and string theory makes predictions about enumerative invariants of
moduli spaces of maps from these and their degenerations to $X$, with
boundary mapping to a special Lagrangian submanifold $L\subset X$.
Enumerative predictions are now being produced by expected dualities,
Ooguri--Vafa~\cite{OV}, Labastida--Mari{\~n}o--Vafa~\cite{LMV},
Aganagic--Vafa~\cite{AV} and Aganagic--Klemm--Vafa~\cite{AKV}.  

Our goals in the present paper are twofold.  First, we outline the
mathematical framework for a theory of enumerative geometry of stable
maps from bordered Riemann surfaces with boundary mapping to a
Lagrangian submanifold (or more simply, with Lagrangian boundary
conditions).  In this direction, we give some foundational results and
propose others which we believe to be true.  Second, we support our
hypotheses by carrying out a computation which under our assumptions
completely agrees with the prediction of \cite{OV} for the
contribution of multiple covers of a disc!  

We were led to modify our conceptual framework after learning of the
results of \cite{AKV}, where it is clarified that
the enumerative invariants are not intrinsic to the geometry and necessarily
depend on an additional discrete parameter.  Happily, our modified 
calculations coincide up to sign with the results of \cite{AKV}, resulting 
in a doubly infinite collection of checks.

Our calculations carried out on prestable bordered Riemann surfaces, as well 
as on the prestable complex algebraic curves which arise from doubling the
bordered surfaces.

Here is an overview of the paper.  In \fullref{enummult}, we give
some background from Gromov--Witten theory and outline some desired
properties for a similar theory of stable maps with Lagrangian
boundary conditions.  In \fullref{bm} we describe bordered Riemann
surfaces and doubling constructions that we will use.  In
\fullref{stable} we discuss stable maps with Lagrangian boundary
conditions and the idea of a virtual fundamental class on moduli
spaces of stable maps.  In \fullref{ta} we describe $U(1)$ actions
on these moduli spaces and perform calculations of the weights of
$U(1)$ representations arising in deformations and obstructions of
multiple covers of a disc.  In \fullref{orientation} we give
orientations on the bundles we need, so that an Euler class can be
defined.  The main results, \fullref{contribution} and
\fullref{ovresult}, are formulated in \fullref{mainresults}.
The proofs are completed in \fullref{calculation}.

We hope that our techniques will prove useful in other situations including
Floer homology and homological mirror symmetry.

Other approaches to defining enumerative invariants of stable maps
with Lagrangian boundary conditions are possible.  M.\ Thaddeus has an
interesting approach using Gromov--Witten theory.  This leads to a
computation in agreement with the $g=0$ case of
\fullref{contribution} Thaddeus~\cite{T}. J Li and Y\,S Song use
relative stable morphisms in their approach. The result of their computation
coincides with the result of \fullref{ovresult} \cite{LS}.

\medskip
\textbf{Acknowledgements}\qua
It is a pleasure to thank J Bryan,
T Graber, J Li, W Li, G Liu, K Liu, C McMullen, Y-G Oh,
J Starr, C\,H Taubes, M Thaddeus, C Vafa, X\,W Wang, and S-T Yau
for helpful conversations.
This research has been supported in part by NSF grant DMS-0073657.
The research of the first author is supported in part by NSA grant
MDA904-00-1-052. 

\section{Enumerative invariants and multiple covers}
\label{enummult}

\subsection{Gromov--Witten invariants and multiple covers of $\PP$}
\label{gwmult}

We start by reviewing ideas from Gromov--Witten theory
relating to multiple covers
which will be helpful for fixing ideas.

Let $X$ be a Calabi--Yau 3--fold and $C\subset X$ a smooth
rational curve with normal 
  bundle $N=N_{C/X}\cong \OO$. Set $\beta=[C]\in H_2(X;\Z)$.
  The virtual dimension of $\MCY$ is $0$. Fixing 
an isomorphism $i\co \PP\simeq C$,
there is a natural embedding $j\co \MPP\rightarrow \MCY$ given by $j(f)=
i\circ f$, whose image is a connected component of $\MCY$, which  we denote
by $\MCY_C$.  The virtual dimension of $\MPP$ is
$2(d+g-1)$.  

Let $\pi\co  \mathcal{U}\rightarrow \MPP$ be the universal family of
stable maps, with $\mu\co \mathcal{U}\rightarrow \PP$ the evaluation 
map. 
Then \[ j_*\left(c(R^1\pi_*\mu^*N\right)\cap
[\MPP]^{\mathrm{vir}})=[\MCY_C]^{\mathrm{vir}} \] where $c(R^1\pi_*\mu^* N)$ is the
top Chern operator of the rank $2(d+g-1)$ obstruction bundle $R^1\pi_*\mu^*
N$ over $\MPP$ and $[\MCY_C]^{\mathrm{vir}}$ denotes the restriction of 
$[\MCY]^{\mathrm{vir}}$ to $[\MCY_C]$.  The contribution $C(g,d)$  
to the genus $g$ Gromov--Witten invariant $N^g_{d[C]}$ of $X$ from
degree $d$ multiple covers of $C$ is defined to be the degree of
$[\MCY_C]^{\mathrm{vir}}$.  It follows that
\[
C(g,d)=\int_{\left[\MPP\right]^{\mathrm{vir}}}c(R^1\pi_*\mu^* N).
\]
Alternatively, 
$C(g,d)$ can be viewed as the Gromov--Witten invariants of the total space of 
$N_{C/X}$, which is a non-compact Calabi--Yau 3--fold.  This follows because
an analytic neighborhood of $C\subset X$ is isomorphic to an analytic
neighborhood of the zero section in $N$, while 
stable maps to $C$ cannot deform off $C$ inside $X$.

The truncated genus $g$ prepotential is 
\begin{equation}
\label{gwpotg}
F_g(t)=\sum_{d=1}^\infty C(g,d)e^{-dt}
\end{equation}
and the all genus truncated potential is 
\begin{equation}
\label{gwpot}
F(\lambda,t)=\sum_{g=0}^\infty \lambda^{-\chi(\Si_g)} F_g(t)
     =\sum_{g=0}^\infty \lambda^{2g-2} F_g(t),
\end{equation}
where $\chi(\Si_g)=2-2g$ is the Euler characteristic of a smooth
Riemann surface of genus $g$.\footnote{The appearance of the Euler
characteristic in the exponent is familiar from considerations of the
topological string, where $\lambda$ is the string coupling constant.}
These potentials are ``truncated'' because they do not contain the
contribution from constant maps.  \ref{gwpot} is viewed as the
contribution of $C$ to the full potential
\begin{equation}
\label{fullpot}
\sum_{g,\beta} \lambda^{2g-2}N^g_\beta q^\beta.
\end{equation}
In \ref{fullpot}, $q^\beta$ is a formal symbol satisfying $q^\beta
q^{\beta'}= q^{\beta+\beta'}$ for classes $\beta,\beta'\in H_2(X,\Z)$.
Alternatively, under suitable convergence hypotheses
$q^\beta$ can be identified with\break $\mathrm{exp}
(-\int_\beta\omega)$, where $\omega$ is the K\"ahler form of $X$,
as discussed in Cox--Katz
\cite{CK}.\footnote{A different constant is used here in the exponent since
one of our goals is to compare to the result of Ooguri--Vafa~\cite{OV}.}

Using the localization formula for the virtual fundamental classes
proven in Graber--Pandharipande~\cite{GP}, it is shown in
Faber--Pandharipande~\cite{FP} that
\[ C(g,d)=d^{2g-3}\sum_{\begin{array}{c}g_1+g_2=g\\g_1,g_2\geq 0
                          \end{array}}b_{g_1}b_{g_2} \]
where $d>0$, and
\[  b_g=\left\{ \begin{array}{ll} 
                  1,&g=0\\
                 \int_{\overM_{g,1}} \psi_1^{2g-2}\lambda_g,&g>0 
                 \end{array}. \right. \]
Therefore
  \begin{eqnarray*}
   F(\lambda,t)&=&\sum_{g=0}^\infty \lambda^{2g-2} \sum_{d=1}^\infty d^{2g-3}
                \sum_{\begin{array}{c}g_1+g_2=g\\g_1,g_2\geq 0
                      \end{array}} b_{g_1} b_{g_2} e^{-dt}\\
   &=& \sum_{d=1}^\infty \frac{e^{-dt}}{\lambda^2 d^3}\sum_{g=0}^\infty 
       \left( \sum_{\begin{array}{c}g_1+g_2=g\\g_1,g_2\geq 0
                    \end{array}} b_{g_1} b_{g_2} \right) (\lambda d)^{2g}\\
   &=& \sum_{d=1}^\infty\frac{e^{-dt}}{\lambda^2 d^3}
       \left( \sum_{g=0}^\infty b_g(\lambda d)^{2g}\right)^2.
  \end{eqnarray*}
It is also shown in \cite{FP} that  
   \[ \sum_{g=0}^\infty b_g u^{2g}=\frac{u/2}{\sin(u/2)}. \]
Therefore 
\[ F(\lambda,t)= \sum_{d=1}^\infty \frac{e^{-dt}}
      {d\left( 2\sin(\lambda d/2)\right)^2} \]
This is the multiple cover formula for the sphere predicted in 
Gopakumar--Vafa~\cite{gov}.

\subsection{Vafa's enumerative invariants.}
\label{vafa}

Several works of Vafa and collaborators deal with enumerative
invariants arising from open string theory, and compute these
invariants in several cases \cite{OV,LMV,AV,AKV}.  It is of fundamental
importance to formulate these invariants as precise mathematical
objects.  We do not do that here, but instead describe properties of
the formulation, much in the same spirit as was done for Gromov--Witten
theory in Kontsevich~\cite{kont}.  These properties will be sufficient to verify
the multiple cover formula of \cite{OV}, to be described at the end of this
section.  The proofs of these properties are left for later work.

Consider a Calabi--Yau 3--fold $X$, and let $L\subset X$ be a special
Lagrangian submanifold.  Fix non-negative integers $g,h$ and a relative
homology class $\beta\in H_2(X,L,\Z)$.  Choose classes $\gamma_1,\ldots,
\gamma_h\in H_1(L,\Z)$ such that $\sum_i\gamma_i=\partial\beta$, where
$$\partial\co H_2(X,L,\Z)\to H_1(L,\Z)$$ 
is the natural map.  Then Vafa's
enumerative invariants are rational numbers $N^{g,h}_{\beta;\gamma_1,\ldots,
\gamma_h}$ which roughly speaking count
continuous maps $f\co  (\Si, \bS)\rightarrow (X, L)$ where
\begin{itemize}
\item $(\Si, \bS)$ is a bordered Riemann surface of genus $g$, whose 
boundary $\bS$ consists of $h$ oriented circles $R_1,\ldots,R_h$.%
\footnote{More precisely, the bordered Riemann surfaces can be
prestable, and the maps $f$ stable.  See Sections~\ref{stabbord} and
\ref{stable}.}
\item $f$ is holomorphic in the interior of $\Si$.
\item $f_*[\Si]=\beta$.
 \item $f_*[R_i]=\gamma_i$.
\end{itemize}
\noindent

To define such an invariant, we will need several ingredients:
\begin{itemize}
\item A compact moduli space 
$\MX$ of stable maps which compactify the space of maps just described.  
\item An orientation on this moduli space.
\item A virtual fundamental class $[\MX]^\mathrm{vir}$ of dimension~0.
\end{itemize}
Then we simply put
\[
N^{g,h}_{\beta;\gamma_1,\ldots,\gamma_h} = \mathrm{deg}[\MX]^\mathrm{vir}.
\]
In \cite{AKV}, it has come to light that these invariants are not completely
intrinsic to the geometry but depend on an additional discrete parameter.
In the context of the dual Chern--Simons theory, a framing is required to 
obtain a topological quantum field theory Witten~\cite{Witten}, and the
ambiguity we refer to arises precisely from this choice.  In the
context of enumerative geometry, this translates into additional
structure required at the boundary of $\MX$, an extension in a sense
of the additional data to be imposed on a complex bundle on a bordered
Riemann surface in order to define its generalized Maslov index (see 
\fullref{tcxdouble}).  In our application, this data can be deduced
from the choice of a normal vector field to $L$.

We do not attempt to describe this circle of ideas here in detail, leaving
this for future work.  Instead, we content ourselves with showing how
different choices lead to different computational results.  We came
to appreciate this point after the authors of \cite{AKV} urged us to look
for an ambiguity in the geometry, and T Graber 
suggested to us that our computations in an earlier draft might depend on
choices of a certain torus action that we will describe later.  An analogous
point was also made in Li--Song~\cite{LS}.  We will explain the relationship between
the choice of torus action and the additional geometric data in 
\fullref{outline}.

We introduce a symbol $q^\beta$ for each $\beta\in H_2(X,L,\Z)$ satisfying
$q^\beta q^{\beta'}=q^{\beta+\beta'}$.  For each $\gamma\in H_1(R_i,\Z)$
we introduce a symbol $q_i^\gamma$ satisfying $q_i^{\gamma}q_i^{\gamma'}=
q_i^{\gamma+\gamma'}$.\footnote{Choosing a generator $g_i\in H_1(R_i,\Z)$
corresponding to the orientation of $R_i$, we can instead introduce
ordinary variables $y_1,\ldots,y_h$ and write $y_i^{d_i}$ in place of
$q_i^{d_ig_i}$.}

We put
\begin{eqnarray*}
F_{g;h}(q^\beta,q_1,\ldots,q_h)&=&
\sum_{\beta}\sum_{\begin{array}{c}\gamma_1+\ldots+\gamma_h=\partial \beta\\
                           \gamma_1,...,\gamma_h\ne 0\end{array}} 
N^{g,h}_{\beta;\gamma_1,\ldots,\gamma_h}
    q_1^{\gamma_1}\cdots q_h^{\gamma_h} q^{\beta}.\\
\end{eqnarray*}
The potential can then be defined to be
\begin{eqnarray*}  
F(\lambda,q^\beta,q_1,\ldots,q_h)&=&\sum_{g=0}^\infty\sum_{h=1}^\infty
\lambda^{-\chi(\Si_{g;h})}F_{g;h}(q^\beta,q_1,\ldots,q_h)\\
&=&\sum_{g=0}^\infty\sum_{h=1}^\infty
\lambda^{2g+h-2}F_{g;h}(q^\beta,q_1,\ldots,q_h),
\end{eqnarray*}
In the above, $\chi(\Si_{g;h})=2-2g-h$ is the Euler characteristic of a
bordered Riemann surface of genus $g$ with $h$ discs removed.

Actually, the moduli space $\MX$ and virtual fundamental class should be 
defined for more general varieties $(X,J,L)$ where $X$ is a symplectic 
manifold, $J$ is a compatible almost complex structure, and $L$ is Lagrangian.
In general however, the virtual dimension will not be 0, as will be
discussed in \fullref{slag}, illustrated for a disc.  The
orientation will be described in \fullref{orientation}, and the
virtual fundamental class will be discussed in Sections~\ref{slag} and
\ref{calculation}.

In many practical examples, $N^{g,h}_{\beta;\gamma_1,\ldots,\gamma_h}$
should be computable using a torus action and a modification of
virtual localization introduced in \cite{GP}.  This is in fact how we
do our main computation, \fullref{contribution}.

In our situation, the computation of the virtual fundamental class is
simplified by exhibiting it as the Euler class of a computable
obstruction bundle.  The Euler class will be computed in
\fullref{calculation}.

\section{Bordered Riemann surfaces and their moduli spaces}
\label{bm}
  
In this section, we describe bordered Riemann surfaces and their
moduli spaces.  We also describe the double of a bordered Riemann
surface and related doubling constructions.  These notions are
connected by the idea of a symmetric Riemann surface, which is a
Riemann surface together with an antiholomorphic involution.  The
considerations developed in this section play a fundamental role in
later computations, as we will frequently alternate between the
bordered Riemann surface and symmetric Riemann surface viewpoints.  For
the convenience of the reader, we have included proofs of well known
results when an appropriate reference does not exist, as well as
related foundational material.

 \subsection{Bordered Riemann surfaces}   
 \label{border}

  This section is a modification of Alling and Greenleaf
\cite[Chapter 1]{AG}.   

  \begin{df}
   Let $A$ and $B$ be nonempty subsets of 
   $\C^+=\{z\in \C\mid\mathrm{Im}z\geq 0 \}$.
   A continuous function $f\co A\rightarrow B$ is \emph{holomorphic} 
   on $A$ if it extends to a holomorphic function $\tilde{f}\co  U\rightarrow\C$,
   where $U$ is an open neighborhood of $A$ in $\C$.
  \end{df}
   
  \begin{thm}[Schwartz reflection principle] 
\begin{sloppypar}
   Let $A$ and $B$ be nonempty subsets of 
   $\C^+=\{z\in \C\mid\mathrm{Im}z\geq0\}$.
   A continuous function $f\co A\rightarrow B$ is holomorphic if it is
   holomorphic on the interior of $A$ and satisfies 
   $f(A\cap\R)\subset B\cap\R$.
\end{sloppypar}
  \end{thm}
   
  \begin{df}
   A \emph{surface} is a Hausdorff, connected, topological space $\Si$
   together with a family $\mathcal{A}=\{(U_i, \phi_i)\mid i\in I\}$ 
   such that $\{U_i\mid i\in I\}$ is an open covering of $\Si$ and each
   map $\phi_i\co U_i\rightarrow A_i$ is a homeomorphism onto an open 
   subset $A_i$ of  $\C^+$. $\mathcal{A}$ is called a 
   \emph{topological atlas} on $\Si$, and each pair $(U_i, \phi_i)$ 
   is called a \emph{chart} of $\mathcal{A}$. The boundary of $\Si$
   is the set
\[ 
   \bS=\{x\in \Si\mid\exists\,i\in I \textup{ s.t. } x\in U_i,
       \phi_i(x)\in\R, \phi_i(U_i)\subseteq \C^+ \}.
\]
   The mappings $\phi_{ij}\equiv \phi_i\circ\phi_j^{-1}\co \phi_j(U_i\cap
   U_j)\rightarrow \phi_i(U_i\cap U_j) $ are surjective
   homeomorphisms, called the \emph{transition functions} of
   $\mathcal{A}$. The atlas $\mathcal{A}$ is called a
   \emph{holomorphic atlas} if all its transition functions are
   holomorphic.  
  \end{df}
  
  \begin{df}
   A \emph{bordered Riemann surface} is a compact surface with nonempty
   boundary equipped with the holomorphic structure induced by a 
   holomorphic atlas on it. 
  \end{df}
 
  \begin{rem}\label{boundary} 
   A Riemann surface is canonically oriented by the holomorphic (complex)
   structure.
   In the rest of this paper, the boundary circles $R_i$ of a bordered
   Riemann surface $\Si$ with boundary $\bS=R_1\cup\ldots\cup R_h$
   will always be given the orientation induced by the
   complex structure, which is a choice of
   tangent vector to $R_i$ such that the basis
   (the tangent vector of $R_i$, inner normal) for the real tangent space
   is consistent with the orientation of $\Si$ induced by the complex 
   structure.
  \end{rem}

  \begin{df}
   A \emph{morphism} between bordered Riemann surfaces $\Si$ and $\Si'$ is
   a continuous map $f\co (\Si, \bS)\rightarrow (\Si', \bS')$
   such that for any $x\in\Si$ there exist analytic charts $(U,\phi)$ and 
   $(V,\psi)$ about $x$ and $f(x)$ respectively, and an analytic function 
   $F\co \phi(U) \rightarrow \C$ such that the following diagram commutes:
\[ 
   \begin{CD}
    U@>f>>V\\@V{\phi}VV @VV{\psi}V\\ \phi(U)@>F>>\C 
   \end{CD}
\] 
  \end{df} 

  A bordered Riemann surface is topologically a sphere with $g\geq 0$ handles 
  and with $h>0$ discs removed. Such a bordered Riemann surface is said to
  be of type $(g;h)$.

 \subsection{Symmetric Riemann surfaces}
 \label{symmetric}
  
  \begin{df} 
   A \emph{symmetric Riemann surface} is a Riemann surface $\Si$ together
   with an antiholomorphic involution $\si\co \Si\rightarrow \Si$. 
   The involution $\si$ is called the \emph{symmetry} of $\Si$. 
  \end{df}
  
  \begin{df}\label{symmor} 
   A \emph{morphism} between symmetric Riemann surfaces $(\Si,\si)$ and
   $(\Si', \si')$ is a holomorphic map $f\co \Si\rightarrow\Si'$ such that 
   $f\circ\si=\si'\circ f$.
  \end{df}

  A compact symmetric Riemann surface is topologically a compact orientable
  surface without boundary $\Si$ together with an orientation reversing 
  involution $\si$, which is classified by the following three invariants:
  \begin{enumerate}      
   \item The genus $\tilde{g}$ of $\Si$.
   \item The number $h=h(\si)$ of connected components of $\Si^\si$, the
         fixed locus of $\si$.
   \item The index of orientability, $k=k(\si)\equiv 2-$the number of
         connected components of $\Si \backslash \Si^\si$.
  \end{enumerate}
  
  These invariants satisfy:
  \begin{itemize}
   \item $0\leq h\leq \tilde{g}+1$.
   \item For $k=0$, we have $h>0$ and $h\equiv \tilde{g}+1$ (mod $2$).
   \item For $k=1$, we have $0\leq h \leq \tilde{g}$.
  \end{itemize}

  The above classification was realized already by Felix Klein (see
  eg, Klein~\cite{Klein}, Weichhold~\cite{W},
  Sepp{\"a}l{\"a}~\cite{S}). This classification is probably better
  understood in terms of the quotient $Q(\Si)=\Si/\bra\si\ket$, where
  $\bra\si\ket=\{id, \si\}$ is the group generated by $\si$.  The
  quotient $Q(\Si)$ is orientable if $k=0$ and nonorientable if $k=1$,
  hence the name ``index of orientability''.  Furthermore, $h$ is the
  number of connected components of the boundary of $Q(\Si)$.  If
  $Q(\Si)$ is orientable, then it is topologically a sphere with
  $g\geq 0$ handles and with $h>0$ discs removed, and the invariants
  of $(\Si,\si)$ are $(\tilde{g},h,k)=(2g+h-1,h,0)$.  If $Q(\Si)$ is
  nonorientable, then it is topologically a sphere with $g>0$
  crosscaps and with $h\geq 0$ discs removed, and the invariants of
  $\Si$ are $(\tilde{g},h,k)=(g+h-1,h,1)$.
 
  From the above classification we see that symmetric Riemann surfaces of a
  given genus $\tilde{g}$ fall into $[\frac{3\tilde{g}+4}{2}]$ topological types.
   
 \subsection{Doubling constructions}
 \label{double}
   
  \subsubsection{The complex double of a bordered Riemann surface}
  \label{doublers}
    
   \begin{thm}\label{doublesurface}
    Let $\Si$ be a bordered Riemann surface. There exists a double cover
    $\pi\co  \Si_\C\rightarrow \Si$ of $\Si$ by a compact Riemann surface
    $\Si_\C$ and an antiholomorphic involution $\si\co \Si_\C\rightarrow\Si_\C$ 
    such that $\pi\circ \si=\pi$. There is a holomorphic embedding 
    $i\co \Si\rightarrow \Si_\C$ such that $\pi\circ i$ is the identity map.
    The triple $(\Si_\C,\pi,\si)$ is unique up to isomorphism.
   \end{thm}
\proof
    See \cite{AG} for the construction of $\Si_\C$. The rest of
    \fullref{doublesurface} is clear from the construction.\endproof

   \begin{df}
    We call the triple $(\Si_\C,\pi,\si)$ in \fullref{doublesurface} the
    \emph{complex double} of $\Si$, and $\overSi=\sigma(i(\Sigma))$
    the \emph{complex conjugate} of $\Si$.
   \end{df}
  
   If $\Si$ is a bordered Riemann surface of type $(g;h)$, then
   $(\Si_\C,\si)$ is a compact symmetric Riemann surface of type
   $(2g+h-1,h,0)$, and $\Si=Q(\Si_\C)$. The two connected components
   of $\Si_\C\backslash(\Si_\C)^{\si}$ are $i(\inS)$ and
   ${\overSi}^\circ=\si\circ i(\inS)$, where $\inS$ denotes the
   interior of $\Si$ and $(\Si_\C)^{\si}$ is the fixed locus of $\si$.

   A bordered Riemann surface of type $(0;1)$ is the disc, and its complex
   double is the projective line; a bordered Riemann surface of type $(0;2)$
   is an annulus, and its complex double is a torus. In the above two cases,
   $\Si$ and $\overSi$ are isomorphic as bordered Riemann surfaces,   
   which is in general not true for bordered Riemann surfaces of other
   topological types.

  \subsubsection{The complex double of a map}
  \label{cxdm}

   \begin{thm}\label{doublemap}
    Let $X$ be a complex manifold with an antiholomorphic involution
    $A\co  X\rightarrow X$, and $L=X^A$ be the fixed locus of $A$.
    Let $\Si$ be a bordered Riemann surface, and 
    $f\co (\Si,\bS)\rightarrow (X,L)$ be a continuous map which is holomorphic
    in the interior of $\Si$.  We identify $\Si$ with its image under $i$ in
    $\Si_\C$. Then $f$ extends to a holomorphic map 
    $f_\C\co \Si_\C\rightarrow X$ such that $A\circ f_\C=f_\C\circ \si$.
   \end{thm}
\proof
    For any $p\in \Si_\C$ define
\[
    f_\C(p)=\left\{ \begin{array}{ll}f(p)&\textup{ if } p\in \Si\\
                   A\circ f\circ \si(p)&\textup{ if } p\notin \Si
                  \end{array} \right.
\]
    Then $f_\C$ is holomorphic by the Schwartz reflection principle, and
    the rest of \fullref{doublemap} is clear. \endproof

   \begin{df}
    We call the holomorphic map $f_\C\co \Si_\C\rightarrow X$ in 
    \fullref{doublemap}
    the \emph{complex double} of $f\co (\Si,\bS)\rightarrow (X,L)$.
   \end{df}
    
   \begin{rem}
\label{notdouble}
    Let $X$ be a complex manifold with an antiholomorphic involution
    $A\co  X\rightarrow X$, and let $L=X^A$ be the fixed locus of $A$.
    Let $\Si$ be a bordered Riemann surface, and let 
    $g\co \Si_\C\rightarrow X$ be a holomorphic map such that 
    $A\circ g=g\circ \si$. It is \emph{not} always true that $g=f_\C$ for 
    some continuous map $f\co (\Si,\bS)\rightarrow (X,L)$ which is holomorphic
    on $\Si^\circ$. For example, let $X=\PP$ and $A(z)={1}/{\bar{z}}$,
    where $z$ is the affine coordinate of $\C\subset \PP$. 
    Then $L=\{|z|=1\}\cong S^1$. Let
    $\Si=D^2=\{|z|\leq 1\}\subset \Si_\C=\PP$ and put
    $\si(z)={1}/{\bar{z}}$. A degree $2$ holomorphic map 
    $g\co \PP\rightarrow \PP$ satisfies $A\circ g=g\circ \si$ if the set
    of two branch points is invariant under $A$, and it is the complex double
    of some continuous map $f\co (D^2,\partial D^2)\rightarrow (\PP, S^1)$
    which is holomorphic on $\Si^\circ$ if the set of the two branch points 
    of $g$ is of the form $\{z_0, {1}/{\bar{z}_0}\}$. In particular,  
    if the branch points of $g$ are two distinct points on $S^1$, then 
    $A\circ g=g\circ\si$, but $g$ is not the complex double of some
    continuous map $f\co (D^2,\partial D^2)\rightarrow (\PP, S^1)$
    which is holomorphic on $\Si^\circ$.
   \end{rem}
     
  \subsubsection{The topological complex double of a complex vector bundle}
\label{tcxdouble}

   We start by reviewing some facts about totally real subspaces of $\C^n$, 
   following Oh~\cite{O}. 

   A real subspace $V$ of $\C^n$ is \emph{totally real} (w.r.t the standard
   complex structure of $\C^n$)  if $\dim_\R V=n$ and $V\cap iV=\{0\}$.
   Define
\[ 
   \mathcal{R}_n\equiv\{V\mid\textup{V is a totally real subspace of } \C^n\}.
\] 
   $GL(n,\C)$ acts transitively on $\mathcal{R}_n$, and the isotropy group of 
   $\R^n\subset \C^n$ is $GL(n,\R)$, so 
   $\mathcal{R}_n\cong GL(n,\C)/GL(n,\R)$. Concretely, this means that
   any totally real subspace 
   $V\subset \C^n$ is of the form $A\cdot\R^n$ for some $A\in GL(n,\C)$, and
   $A_1\cdot\R^n=A_2\cdot \R^n$ if and only if $A_2^{-1}A_1\in GL(n,\R)$, or 
   equivalently, $A_1\bar{A}_1^{-1}=A_2\bar{A}_2^{-1}$.  By
   \cite[Proposition 4.4]{O}, the map
   \begin{eqnarray*}
    B\co \mathcal{R}_n\cong GL(n,\C)/GL(n,\R)&\longrightarrow&
    \widetilde{\mathcal{R}}_n\equiv\{D\in GL(n,\C)\mid D\bar{D}=I_n\}\\
    A\cdot \R^n&\mapsto& A\bar{A}^{-1} 
   \end{eqnarray*}
   is a diffeomorphism, where $I_n$ is the $n\times n$ identity matrix. 
   The \emph{(generalized) Maslov index} $\mu(\gamma)$ of an oriented loop 
   $\gamma\co S^1\rightarrow \mathcal{R}_n$ is defined to be the degree of
   the map $\phi=\det\circ B\circ \gamma\co  S^1\rightarrow U(1)$, where 
   $U(1)=\{e^{i\theta}\mid\theta\in \R\}$ is oriented by 
   ${\partial}/{\partial\theta}$. 

   The (generalized) Maslov index gives an explicit way to detect homotopy
   of loops in $\mathcal{R}_n$; in fact $\pi_1(\mathcal{R}_n)\cong \Z$, and
   two loops $\gamma_1, \gamma_2\co S^1\rightarrow \mathcal{R}_n$ are homotopic
   if and only if $\mu(\gamma_1)=\mu(\gamma_2)$.

   A real subspace $V$ of $\C^n$ is Lagrangian with respect to
   the standard symplectic structure on $\C^n$ if $\dim_\R V=n$ and
   $V$ is orthogonal to $iV$ with respect to the standard inner product on
   $\R^{2n}\cong \C^n$.  The space $\mathcal{L}_n$ of Lagrangian
   subspaces of $\C^n$ is a submanifold of $\mathcal{R}_n$, and its
   image under $B$ is the submanifold
\[
   \widetilde{\mathcal{L}}_n\equiv\{D\in GL(n,\C)\mid D\bar{D}=I_n,\ D=D^t\}
\]
   of $\widetilde{\mathcal{R}}_n$. The inclusion
   $\mathcal{L}_n\subset\mathcal{R}_n$ is a homotopy equivalence, and
   the restriction of the (generalized) Maslov index to loops in
   $\mathcal{L}_n$ is the usual Maslov index in symplectic geometry
   (McDuff--Salamon \cite[Chapter 2]{MS}).
  
   Let $E$ be a complex vector bundle of rank $n$ over a bordered Riemann
   surface $\Si$, and let $E_\R$ be a totally real subbundle of
   $E|_{\bS}$, so that
   \[ E|_{\bS}\cong E_\R\otimes_\R \C. \]
   We may assume that $h>0$, so
   that $E$ is a topologically trivial complex vector
   bundle.  Fix a trivialization $\Phi\co E\cong \Si\times\C^n$, and let
   $R_1, \ldots, R_h$ be the connected components of $\bS$, with orientation
   induced by the orientation of $\Si$, as explained in 
   \fullref{boundary}. Then $\Phi(E_\R|_{R_i})$ gives
   rise to a loop $\gamma_i\co S^1\rightarrow\mathcal{R}_n$. Setting  
   $\mu(\Phi, R_i)=\mu(\gamma_i)$, we have

  \begin{pro}
  \label{Maslov}
   Let $E$ be a complex vector bundle of rank $n$ over a bordered Riemann
   surface $\Si$, and let $E_\R$ be a totally real subbundle of $E|_{\bS}$.
   Let $R_1, ..., R_h$ be the connected components of $\bS$, with 
   orientation induced by the orientation of $\Si$. Then
   $\sum_{i=1}^h \mu(\Phi,R_i)$ is independent of the choice of 
   trivialization $\Phi\co E\cong \Si\times \C^n$. 
  \end{pro}
\proof
   Let $\Phi_1$, $\Phi_2$ be two trivializations of $E$. Then
   $\Phi_2\circ\Phi_1^{-1}\co \Si\times\C^n \rightarrow\Si\times\C^n$ is
   given by
   $(x,v)\mapsto(x,g(x)v)$ for an appropriate 
   $g\co \Si\rightarrow GL(n,\C)$. Letting
\[ 
   \phi=\det(g\bar{g}^{-1})\co \Si\rightarrow U(1),
\]  
   it follows that
\[ 
   \mu(\Phi_2,R_i)-\mu(\Phi_1, R_i)=\deg(\phi|_{R_i}).
\] 
   Furthermore, $\phi_*[R_i]=\deg(\phi|_{R_i})[U(1)]$, where
   $\phi_*\co H_1(\Si;\Z)\rightarrow H_1(U(1);\Z)$ and $[U(1)]\in
   H_1(U(1),\Z)$ is the fundamental class of $U(1)$.
   Since $[R_1]+\cdots+[R_h]=[\bS]=0\in H_1(\Si;\Z)$, it follows that
   $\sum_{i=1}^h\deg(\phi|_{R_i})=0$, which implies\break
   $\sum_{i=1}^h \mu(\Phi_1, R_i)=\sum_{i=1}^h \mu(\Phi_2, R_i)$, 
   as desired. \endproof 
     
  \begin{df}
   The \emph{Maslov index} of $(E,E_\R)$ is defined by 
\[ 
   \mu(E,E_\R)=\sum_{i=1}^h\mu(\Phi,R_i) 
\]
   where $\Phi\co E\rightarrow \Si\times \C^n$ is any trivialization.
  \end{df}

  \fullref{Maslov} says that the Maslov index is well defined.

\smallskip
  The following theorem is well-known. We include the proof for completeness.
  \begin{thm}
  \label{degree}
   Let $E$ be a complex vector bundle over a bordered Riemann surface
   $\Si$, and let $E_\R$ be a totally real subbundle of $E|_{\bS}$.
   Then there is a complex vector bundle $E_\C$ on $\Si_\C$ together
   with a conjugate linear involution $\tsi\co E_\C\rightarrow E_\C$
   covering the antiholomorphic involution $\si\co \Si_\C\rightarrow\Si_\C$
   such that $E_\C|_\Si=E$ (where $\Si$ is identified with its image 
   under $i$ in $\Si_\C$) and the fixed locus of $\tsi$
   is $E_\R\rightarrow \bS$. Moreover, we have
\[
   \deg E_\C=\mu(E,E_\R).
\]
  \end{thm}
   Here, $\deg E_\C$ means $\deg(c_1(E_\C)\cap[\Si_\C])$ as usual.

\proof
   Let $R_1,\ldots,R_h$ be the connected components
   of $\bS$, and let $N_i\cong R_i\times [0,1)$ be a neighborhood of
   $R_i$ in $\Si$ such that $N_1,\ldots,N_h$ are disjoint. Then
   $(N_i)_\C=N_i\cup \overline{N}_i$ is a tubular neighborhood of $R_i$ in 
   $\Si_\C$, and $N\equiv\cup_{i=1}^h (N_i)_\C$ is a tubular 
   neighborhood of $\bS$ in $\Si_\C$. Let $U_1=\Si\cup N$, 
   $U_2=\overSi\cup N$, so that $U_1\cup U_2=\Si_\C$ and 
   $U_1\cap U_2=N$.
 
   Fix a trivialization $\Phi\co E\cong\Si\times\C^n$, where $n$ is 
   the rank of E. Then $\Phi(E_\R|_{R_i})$ gives rise to a loop
   $B_i\co R_i\rightarrow \widetilde{\mathcal{R}}_n\subset GL(n,\C)$.
   To construct $E_\C\rightarrow \Si_\C$, we glue trivial
   bundles $U_1\times \C^n\rightarrow U_1$ and 
   $U_2\times \C^n\rightarrow U_2$ along $N$ by identifying
   $(x,u)\in (N_i)_\C\times \C^n\subset U_1\times \C^n$
   with $(x,B_i^{-1}\circ p_i(x)u)\in(N_i)_\C\times 
   \C^n\subset U_2\times \C^n$,
   where $p_i\co (N_i)_\C\cong R_i\times(-1,1)\rightarrow R_i$ is the
   projection to the first factor and 
   $B_i^{-1}\co R_i\to\widetilde{\mathcal{R}}_n$ denotes the map $B_i^{-1}(x)=
   (B_i(x))^{-1}$. There is a conjugate linear 
   involution $\tsi\co E_\C\rightarrow E_\C$ given by
   $(x,u)\in U_1\times\C^n\mapsto (\si(x),\bar{u})\in U_2\times\C^n$
   and     
   $(y,v)\in U_2\times\C^n\mapsto (\si(y),\bar{v})\in U_1\times\C^n$.
   It is clear from the above construction that
    $\tsi\co E_\C\rightarrow E_\C$ covers the antiholomorphic 
   involution $\si\co \Si_\C\rightarrow\Si_\C$, and the fixed locus of
   $\tsi$ is $E_\R\rightarrow \bS$.

   By considering the determinant line bundles, we only need to 
   show $\deg E_\C=\mu(E,E_\R)$ for a line bundle $E$. 
   In this case, 
   $B_i\co R_i\rightarrow \widetilde{\mathcal{R}}_1=U(1)\subset GL(1,\C)$.
   Let $F\in \Omega^2(\Si_\C,i\R)$ be the curvature of some 
   $U(1)$ connection on $E$. The proof of~\cite[Theorem 2.70]{MS}
   shows that 
\[ 
   \mu(E,E_\R)=\frac{i}{2\pi}\int_{\Si_\C}F= \deg(E_\C).\eqno{\qed}
\]

\begin{df}
   We call the bundle $E_\C$ as in \fullref{degree} the
   \emph{topological complex double} of $(E,E_\R)\rightarrow (\Si,\bS)$.
  \end{df}

  \subsubsection{The holomorphic complex double of a Riemann--Hilbert bundle}
  \label{cxdouble}
  
  \begin{df}
   Let $A$ be a nonempty subset of $\C^+$.
   A continuous function $f\co A\rightarrow \C^n$ is \emph{holomorphic} on $A$ 
   if it extends to a holomorphic function $\tilde{f}\co  U\rightarrow\C^n$,
   where $U$ is an open neighborhood of $A$ in $\C$.
  \end{df}
  
  \begin{df}
   Let $\Si$ be a bordered Riemann surface, and let $X$ be a complex
   manifold. We say $f\co \Si\rightarrow X$ is \emph{holomorphic} if
   for any $x\in\Si$ there exist holomorphic charts $(U,\phi)$ and 
   $(V,\psi)$ about $x$ and $f(x)$ respectively, and a holomorphic function 
   $F\co \phi(U) \rightarrow \C^n$ such that the following diagram commutes:
\[
   \begin{CD}
    U@>f>>V\\@V{\phi}VV @VV{\psi}V\\ \phi(U)@>F>>\C^n 
   \end{CD}
\] 
  \end{df} 
  
  \begin{df}
   A complex vector bundle $E$ over a bordered Riemann surface $\Si$ is
   \emph{holomorphic} if there is an open covering $\{U_i\mid i\in I\}$ of
   $\Si$ together with trivializations $\Phi_i\co E|_{U_i}\cong U_i\times \C^n$ 
   such that if $U_i\cap U_j\neq\phi$ then
   \begin{eqnarray*} 
   \Phi_{ij}\equiv \Phi_i\circ \Phi_j^{-1}\co  
    (U_i\cap U_j)\times \C^n&\longrightarrow&(U_i\cap U_j)\times \C^n\\
    (x,u)&\mapsto&(x, g_{ij}(x)u)
   \end{eqnarray*}
   for some holomorphic map $g_{ij}\co U_i\cap U_j\rightarrow GL(n,\C)$.  
  \end{df}  

   Let $E$ be a smooth complex vector bundle of rank $n$ over a bordered 
   Riemann surface $\Si$ which is holomorphic on $\inS$.  
   More explicitly, there is an open covering $\{U_i\mid i\in I\}$ of
   $\Si$ together with trivializations $\Phi_i\co E|_{U_i}\cong U_i\times \C^n$ 
   such that if $U_i\cap U_j\neq\phi$ then
   \begin{eqnarray*} 
   \Phi_{ij}\equiv \Phi_i\circ \Phi_j^{-1}\co  
    (U_i\cap U_j)\times \C^n&\longrightarrow&(U_i\cap U_j)\times \C^n\\
    (x,u)&\mapsto&(x, g_{ij}(x)u)
   \end{eqnarray*}
   for some smooth map 
   $g_{ij}\co U_i\cap U_j\rightarrow GL(n,\C)$ which is holomorphic
   on $(U_i\cap U_j)\cap\inS$.  
   We call $\{(U_i,\Phi_i)\mid i\in I\}$ a smooth
   trivialization of $(E,E_\R)\rightarrow (\Si,\bS)$ which is 
   holomorphic on $\inS$.
   
  \begin{thm}\label{holbundle}
   Let $E$ be a smooth complex vector bundle of rank $n$ over a bordered
   Riemann surface $\Si$ which is holomorphic on $\inS$, and let $E_\R$ 
   be a smooth totally real subbundle of $E|_{\bS}$. Then there exists
   a smooth trivialization $\{(U_i,\Phi_i)\mid i\in I\}$ of $(E,E_\R)$ 
   which is holomorphic on $\inS$ such that
\[ 
   \Phi_i(E_\R|_{U_i\cap\bS})=(U_i\cap\bS)\times \R^n\subset U_i\times\C^n
\]     
   whenever $U_i\cap\bS\neq \phi$.
   In particular, $\{(U_i,\Phi_i)\mid i\in I\}$ determines a holomorphic
   structure on $E\rightarrow \Si$ such that $E_\R\rightarrow\bS$ is a 
   real analytic subbundle of $E|_{\bS}\rightarrow \bS$.  
  \end{thm}
\proof
   Let \{$(U_i,\Phi'_i)\mid i\in I$\} be a smooth trivialization of $(E,E_\R)$
   which is holomorphic on $\inS$.
   If $U_i\cap\bS$ is empty, set $\Phi_i=\Phi'_i$. If $U_i\cap \bS$ is 
   nonempty, our strategy is to find a smooth map 
   $h\co U_i\rightarrow GL(n,C)$ which is holomorphic on $U_i\cap \inS$ 
   such that
   $H\circ\Phi_i(E_\R|_{U_i\cap\bS})=(U_i\cap\bS)\times \R^n$,
   where $H\co U_i\times \C^n\rightarrow U_i\times \C^n$, 
   $(x,u)\mapsto (x,h(x)u)$, and set $\Phi_i=H\circ\Phi'_i$. 
    
   By refining the open covering, we may assume that there is
   a holomorphic map $\phi_i\co U_i\rightarrow \C^+$ such that
   $(U_i, \phi_i)$ is a holomorphic chart. Then 
   $\phi_i(U_i\cap\bS)\subset \R$. Let $\phi'_i\co U_i\rightarrow D^2$
   be the composition of $\phi_i$ with the isomorphism $\C^+\cong D^2$.
   We may assume that 
   $U\equiv\phi'_i(U_i)\subset\{z\in D^2\mid \mathrm{Im}z>0\}\equiv T$.
   Then $E_\R$ gives rise to a smooth map 
   $B\co U\cap\partial D^2\rightarrow \widetilde{\mathcal{R}}_n$, and the
   proof will
   be completed if we can find a smooth map $A\co U\rightarrow GL(n,\C)$
   which is holomorphic on the interior of $U$ such that  
   $A\bar{A}^{-1}|_{U\cap\partial D^2}=B$, since $h=A^{-1}\circ\phi'_i$ 
   will have the desired property described in the previous paragraph.

   We extend $B$ to a smooth map 
   $\tilde{B}\co \partial D^2\rightarrow\widetilde{\mathcal{R}}_n$.
   By~\cite[Lemma 4.6]{O}, there is a smooth map 
   $\Theta\co D^2\rightarrow GL(n,\C)$
   which is holomorphic on the interior of $D^2$ such that
   $B(z)=\Theta(z)$diag$(z^{\kappa_1}, \ldots,z^{\kappa_n})
   \overline{\Theta(z)}^{-1}$, for some integers $\kappa_1,\ldots,\kappa_n$.
   Let $r(z)$ be a holomorphic branch of $z^{1/2}$ on $T$.
   Set $A(z)=\Theta(z)$diag$(r(z)^{\kappa_1}\ldots r(z)^{\kappa_n})$ for
   $z\in U$. Then $A$ is smooth and is holomorphic on the interior of $U$,
   and $A\bar{A}^{-1}|_{U\cap\partial D^2}=B$, as desired. 
    
   Finally, if $U_i\cap U_j\neq\phi$, then
   \begin{eqnarray*} 
   \Phi_{ij}\equiv \Phi_i\circ \Phi_j^{-1}\co  
    (U_i\cap U_j)\times \C^n&\longrightarrow&(U_i\cap U_j)\times \C^n\\
    (x,u)&\mapsto&(x, g_{ij}(x)u)
   \end{eqnarray*}
   for some smooth map $g_{ij}\co U_i\cap U_j\rightarrow GL(n,\C)$ 
   which is holomorphic on $U_i\cap U_j\cap\inS$, and 
   $g_{ij}(U_i\cap U_j\cap \bS)\subset GL(n,\R)$.  Note that
   $U_i\cap U_j\cap\inS=U_i\cap U_j$ if and only if 
   $U_i\cap U_j\cap \bS$ is empty.
   By the Schwartz reflection principle, $g_{ij}$ can be extended to  
   a holomorphic map $(g_{ij})_\C\co (U_i\cap U_j)_\C \rightarrow GL(n,\C)$,
   and $g_{ij}|_{U_i\cap U_j\cap \bS}$ is real analytic.
   Therefore, $\{(U_i,\Phi_i)\mid i\in I\}$ determines a holomorphic
   structure on $E\rightarrow \Si$ such that $E_\R\rightarrow\bS$ is a 
   real analytic subbundle of $E|_{\bS}\rightarrow \bS$. \endproof

  \begin{df}
   We call $(E,E_\R)\rightarrow (\Si,\bS)$ together with the holomorphic 
   structure on $E\rightarrow \Si$ determined by 
   $\{(U_i,\Phi_i)\mid i\in I\}$ as in \fullref{holbundle} a 
   \emph{Riemann--Hilbert bundle} over $\Si$.  The trivialization
   $\{(U_i,\Phi_i)\mid i\in I\}$
   is called a \emph{Riemann--Hilbert trivialization} of 
   $(E,E_\R)\rightarrow (\Si,\bS)$.     
  \end{df}
 
   The holomorphic maps $(g_{ij})_\C\co (U_i\cap U_j)_\C \rightarrow GL(n,\C)$
   in the proof of \fullref{holbundle} give a holomorphic vector 
   bundle $E_\C\rightarrow\Si_\C$ together with a holomorphic trivialization 
   $\{((U_i)_\C,(\Phi_{ij})_\C)\mid i\in I\}$.
   There is an antiholomorphic involution $\tsi\co E_\C\rightarrow E_\C$
   such\break\eject that $\tsi(E|_{(U_i)_\C})=E|_{(U_i)_\C}$ and
   \begin{eqnarray*}
    \Phi_i\circ\tsi\circ\Phi_i^{-1}\co
    (U_i)_\C\times \C^n&\longrightarrow&(U_i)_\C\times \C^n\\
    (x,u)&\mapsto&(\si(x),\bar{u})
   \end{eqnarray*}
   It is clear from the above construction that
   $E_\C|_{\Si}=E$ (where $\Si$ is identified with its image
   under $i$ in $\Si_\C$), $\tsi$ covers the antiholomorphic
   involution $\si\co \Si_\C\rightarrow\Si_\C$, and the fixed locus
   of $\tsi$ is $E_\R\rightarrow \bS$.

  \begin{df}
  \label{doublebundle}
   Let $(E,E_\R)\rightarrow (\Si,\bS)$ be a Riemann--Hilbert bundle 
   over a bordered Riemann surface $\Si$, and construct 
   $E_\C\rightarrow \Si_\C$ as above. We call $E_\C\rightarrow\Si_\C$
   the \emph{holomorphic complex double} of the Riemann--Hilbert bundle
   $(E,E_\R)\rightarrow (\Si,\bS)$.
  \end{df}

  \begin{rem}
   Let $(E,E_\R)\rightarrow (\Si,\bS)$ be a Riemann--Hilbert bundle
   over a bordered Riemann surface $\Si$. 
   The underlying topological complex vector bundle of the holomorphic
   complex double $E_\C$ of $(E,E_\R)$ is isomorphic (as a topological
   complex vector bundle) to the topological complex double of
   the underlying topological bundles of $(E,E_\R)$. In particular,
   $\deg E_\C=\mu(E,E_\R)$.

  \end{rem}

\subsection{Riemann--Roch theorem for bordered Riemann surfaces}

  Let $(E,E_\R)\rightarrow (\Si,\bS)$ be a Riemann--Hilbert bundle
  over a bordered Riemann surface $\Si$, and let $(\E,\E_\R)$ denote 
  the sheaf of local holomorphic sections of $E$ with boundary values in $E_\R$.
  Let $\E_\C$ denote the sheaf of local holomorphic sections of $E_\C$, 
  the holomorphic complex double of $(E,E_\R)$. Let $\tilde{U}$ be an open 
  subset of $\Si_\C$. Define $\tsi\co \E_\C(\tilde{U})\rightarrow 
  \E_\C(\si(\tilde{U}))$ by $\tsi(s)(x)=\tsi\circ s \circ \si(x)$ for 
  $x\in \si(\tilde{U})$. In particular, if $U$ is an open set in $\Si$, then
  $\tsi\co \E_\C(U_\C)\rightarrow\E_\C(U_\C)$, and the fixed locus
  $\E_\C(U_\C)^{\tsi}=(\E,\E_\R)(U)$ is a totally real subspace of 
  the complex vector space $\E_\C(U_\C)$.
  
  We are interested in the sheaf cohomology groups $H^q(\Si,\bS,E,E_\R)$
  of the sheaf $(\E,\E_\R)$ on $\Si$ . 
  Here the sheaf cohomology functors are the right derived functors of the
  global section functor from the category of sheaves of $(\Oh,\Oh_\R)$ modules
  on $\Si$ to the category of $\R$ modules, where $(\Oh,\Oh_\R)$ is the sheaf
  of local holomorphic functions on $\Si$ with real boundary values.
  Let $\mathcal{A}^0(E,E_\R)$ denote the
  sheaf of local $C^\infty$ sections of $E$ with boundary values in $E_\R$,
  and let $\mathcal{A}^{0,1}(E)$ denote the sheaf of local $C^\infty$
  $E$--valued $(0,1)$ forms. 
  We claim that
  \begin{equation} \label{fine} 
  0\rightarrow (\E,\E_\R)\rightarrow \mathcal{A}^0(E,E_\R)
    \stackrel{\bar{\partial}}{\rightarrow} \mathcal{A}^{0,1}(E)\rightarrow 0
  \end{equation}
  is a fine resolution of $(\E,\E_R)$. The sheaves
$\mathcal{A}^0(E,E_\R)$ and 
  $\mathcal{A}^{0,1}(E)$ are clearly fine. By definition, the kernel of
  $\bar{\partial}$ is $(\E,\E_\R)$. To check the surjectivity of 
  $\bar{\partial}$, it suffices to check the surjectivity of the stalk
  $\bar{\partial}_x$ at each point $x\in\Si$. For $x\in \inS$, this follows
  from the $\bar{\partial}$--Poincar\'{e} lemma. A neighborhood of $x\in \bS$ 
  in $\Si$ can be identified with an open set $U\subset D^2$ such that
  $U\cap\partial D^2\neq\phi$, and by shrinking $U$ we may assume that
  $(E,E_\R)|_U$ gets identified with the trivial bundle
  $(\C^n,\R^n)\rightarrow (U,U\cap\partial D^2)$. 
  Let $F\co U\rightarrow \C^n$ be a $C^\infty$ function. We need to find 
  a neighborhood $V$ of $x$ in $U$ and a $C^\infty$ function
  $\eta\co V\rightarrow\C^n$ such that
\[
  \left\{\begin{array}{ll}\frac{\partial\eta}{\partial\bar{z}}
   =F(z)&\textup{ on }V^\circ\\ 
   \eta(z)\in\R^n&\textup{ on }V\cap\partial D^2.
  \end{array}\right. 
\]
  Let $\rho$ be a $\C^\infty$ cut-off function whose support is contained 
  in $U$ and satisfies $\rho\equiv 1$ on some neighborhood $V$ of $x$. Then
  $G(z)\equiv\rho(z)F(z)$ defines a $C^\infty$ function $G\co D^2\rightarrow \C^n$.
  It suffices to solve 
\[
  \left\{\begin{array}{ll}\frac{\partial\eta}{\partial\bar{z}}
   =G(z)&\textup{ on }(D^2)^\circ\\ 
  \eta(z)\in \R^n&\textup{ on }\partial D^2.
  \end{array}\right. 
\]
  The above system can be solved (see \cite{V}, \cite{O}). We conclude that
  $\bar{\partial}_x$ is surjective for $x\in \bS$.
 
The resolution \ref{fine} is fine, so the sheaf cohomology of 
  $(\E,\E_\R)$ is given by the cohomology of the two-term elliptic complex
  \begin{equation}\label{Dol}
   0\rightarrow A^0(E,E_\R)\stackrel{\bar{\partial}_{(E,E_\R)}}
   {\longrightarrow}A^{0,1}(E)\rightarrow 0,
  \end{equation}
  where $A^0(E,E_\R)$ is the space of global $C^\infty$ sections of $E$
  with boundary values in $E_\R$, and $A^{0,1}(E)$ is the space of global
  $C^\infty$ $E$--valued $(0,1)$ forms. In other words,
  \begin{eqnarray*}
   H^0(\Si,\bS,E,E_\R)&=&\mathrm{Ker}\,\bar{\partial}_{(E,E_\R)}\\
   H^1(\Si,\bS,E,E_\R)&=&\mathrm{Coker}\,\bar{\partial}_{(E,E_\R)}\\ 
   H^q(\Si,\bS,E,E_\R)&=&0\;\;\;\textup{ for } q>1
  \end{eqnarray*}
  To summarize, the sheaf cohomology of $(\E,\E_\R)$ can be identified
  with the \emph{Dolbeault cohomology}, the cohomology of the 
  \emph{twisted Dolbeault complex} \ref{Dol}. 

  Let $\mathcal{A}=\{U_i\mid i\in I\}$ be an acyclic cover of $\Si$ for the
  sheaf $(\E,\E_\R)$ in the sense that 
\[ 
  H^q(U_{i_1}\cap\cdots\cap U_{i_p}, 
      U_{i_1}\cap\cdots\cap U_{i_p}\cap\bS,E,E_\R)=0
  \textup{ for }q>0, i_1,\ldots,i_p\in I.
\]
  Then the sheaf cohomology group $H^q(\Si,\bS,E,E_\R)$ is isomorphic to 
  the \v{C}ech cohomology group $H^q(\mathcal{A},\E,\E_\R)$, for all $q$.
  We may further assume that $\mathcal{A}\equiv\{(U_i)_\C\mid i\in I\}$ is
  an acyclic cover for the sheaf $\E_\C$ on $\Si_\C$ in the sense that
  the sheaf cohomology groups 
\[ 
  H^q((U_{i_1})_\C\cap\cdots\cap (U_{i_p})_\C,E_\C)=0
  \textup{ for }q>0, i_1,\ldots,i_p\in I.
\]
  Then the sheaf cohomology group $H^q(\Si,E_\C)$ of $E_\C$ is isomorphic to 
  the \v{C}ech cohomology group $H^q(\mathcal{A}_\C,\E)$, for all $q$.
  We have seen that for any open subset $U\subset\Si$, there is a conjugate
  linear involution $\tsi\co \E_\C(U_\C)\rightarrow\E_\C(U_\C)$
  whose fixed locus $\E_\C(U_\C)^{\tsi}= (\E,\E_\R)(U)$.  Note that
  $\tsi$ acts on the \v{C}ech complex $C^\bullet(\mathcal{A}_\C,\E_\C)$,
  and the fixed locus $C^\bullet(\mathcal{A}_\C,\E_\C)^{\tsi}$ is the
  \v{C}ech complex $C^\bullet(\mathcal{A},\E,\E_\R)$. 
  So $\tsi$ acts on $H^q(\Si_\C,E_\C)\cong H^q(\mathcal{A}_\C,\E_\C)$, and
  $H^q(\Si,\bS,E,E_\R)\cong H^q(\Si_\C,E_\C)^{\tsi}$. In particular,
\[
  \dim_\R H^q(\Si,\bS,E,E_\R)=\dim_\C H^q(\Si_\C,E_\C).
\]
 
  \begin{df} 
  \label{index}
   Let $(E,E_\R)\rightarrow (\Si,\bS)$ be a Riemann--Hilbert bundle
   over a bordered Riemann surface $\Si$.
   The \emph{index} of $(E, E_\R)$ is the virtual real vector space 
\[ 
   Ind(E,E_\R)\equiv H^0(\Si,\bS,E,E_\R)- H^1(\Si,\bS,E,E_\R).
\]
   The \emph{Euler characteristic} of $(E,E_\R)$ is the virtual dimension 
\[ 
   \chi(E,E_\R)\equiv\dim_\R H^0(\Si,\bS,E,E_\R)-\dim_\R H^1(\Si,\bS,E,E_\R) 
\]
   of the virtual vector space $Ind(E,E_\R)$.
  \end{df}

  Note that $\chi(E,E_\R)$ is the index of the Fredholm operator 
  $\bar{\partial}_{(E,E_\R)}$. Furthermore, $\chi(E,E_\R)=\chi(E_\C)$. 
  We have the following Riemann--Roch theorem for bordered Riemann surfaces.
  \begin{thm} 
  \label{RR}
   Let $(E,E_\R)\rightarrow(\Si,\bS)$ be a Riemann--Hilbert bundle of 
   rank $n$ over a bordered Riemann surface $\Si$ of type $(g;h)$. Then
\[
   \chi(E,E_\R)=\mu(E,E_\R)+n\chi(\Si),
\]
   where $\chi(\Si)=2-2g-h$ is the Euler characteristic of $\Si$. 
  \end{thm}
\proof
   Let $\tilde{g}=2g+h-1$ be the genus of
   $\Si_\C$. Then 
   \begin{eqnarray*}
    \chi(E,E_\R)&=&\chi(E_\C)\\
                &=&\deg E_\C +n(1-\tilde{g}) \\
                &=&\mu(E,E_\R)+ n(2-2g-h) \\
                &=&\mu(E,E_\R)+ n \chi(\Si), 
   \end{eqnarray*}
   where the second equality follows from
   Riemann--Roch and the third equality comes from
   \fullref{degree}. \endproof

The computations in the following examples will be useful in the sequel.
    
\begin{example} \label{double1}
Let $m$ be an integer. Consider the line bundle $(L(m),L(m)_\R)$ over
   $(D^2,S^1)$, where $L(m)$ is the trivial bundle, and the fiber of
   $L(m)_\R$ over $z=e^{i\theta}$ is
   $ie^{(im\theta/2)}\R\subset\C$.  Then
   $\mu(L(m),L(m)_\R)=m$. The complex double $D^2_\C$ can be identified
   with $\PP$. The antiholomorphic
   involution is $\si(z)=1/\bar{z}$ in an affine coordinate $z$ on
   $\PP$, and $L(m)_\C$ can be identified with $\Oh_{\PP}(m)$.  

   Explicitly, let $(z,u)$ and
   $(\tilde{z},\tilde{u})$ be the two charts of the total
   space of $\Oh_{\PP}(m)$, related
   by $(\tilde{z},\tilde{u})=(1/z,-z^{-m}u)$. There is an
   antiholomorphic involution 
   \begin{eqnarray*}
    \tsi\co \Oh_{\PP}(m)&\longrightarrow&\Oh_{\PP}(m)\\
     (z,u)&\longmapsto&(\frac{1}{\bar{z}},-\bar{z}^{-m}\bar{u})
   \end{eqnarray*}
   expressed in terms of the first chart (it clearly antiholomorphic
   over $z=0$ as well, as $\tsi(z,u)=(\bar{z},\bar{u})$ when the image
   is expressed in the coordinates of the second chart). 

   Note that $\tsi$ covers $\si\co \PP\rightarrow \PP$,
   and the fixed locus of $\tsi$ is $L(m)_\R$. 
   For any nonnegative integer $m$ we have 
   \begin{eqnarray*}
    H^0(\PP,\Oh_{\PP}(m))&=&\left\{\left.\sum_{j=0}^m a_j
     X_0^{m-j}X_1^j\;\right|\; a_j\in\C\right\},\\
    && \tsi(a_0, a_1, \ldots, a_m)=(-\bar{a}_m,
     \ldots,-\bar{a}_1,-\bar{a}_0)\\ 
    H^1(\PP,\Oh_{\PP}(m))&=&0\\
    H^0(\PP,\Oh_{\PP}(-m-1))&=&0\\ 
    H^1(\PP,\Oh_{\PP}(-m-1))&=&
    \left\{\left.\sum_{j=1}^m\frac{a_j}{X_0^{m-j}X_1^j}\;\right|\;
     a_j\in\C\right\},\\
    && \tsi(a_1, a_2,\ldots, a_m)=(-\bar{a}_m,\ldots, -\bar{a}_2,-\bar{a}_1),
   \end{eqnarray*}
   where $(X_0, X_1)$ are homogeneous coordinates on $\PP$,
   related to $z$ by $z=\frac{X_1}{X_0}$, and the classes in $H^1$ are
   expressed as \v{C}ech cohomology classes using the standard open cover 
   of $\PP$. Therefore,
   \begin{eqnarray*} 
    H^0(D^2,S^1, L(m), L(m)_\R)&=&\left\{\left.f(z)=\sum_{j=0}^m
     a_j z^j\;\right|\;a_{m-j}=-\bar{a}_j\right\}\\ 
    H^1(D^2,S^1, L(m), L(m)_\R)&=&0\\
    H^0(D^2,S^1, L(-m-1), L(-m-1)_\R)&=&0\\ 
   H^1(D^2,S^1, L(-m-1), L(-m-1)_\R)&\cong&\left\{\!\left.
     f(z)=\sum_{j=1}^m\frac{a_j}{z^j}\;\right|\;
     a_{m+1-j}=-\bar{a}_j\!\right\}.
   \end{eqnarray*} 
  \end{example}

\begin{example} \label{double2}
Consider next the bundle $(N(d),N(d)_\R)$
   on $(D^2,S^1)$, where $N(d)$ is the trivial rank~2 bundle,
   and the fiber of $N(d)_\R$ over $z\in S^1$ is
\[
   \left\{(u,v)\in\C^2\mid v=\bar{z}^d\bar{u}\right\}.
\]
   Then $\mu(N(d), N(d)_\R)=-2d$, and $N(d)_\C=\OOd$.
   Let $(z,u,v)$ and $(\tilde{z}, \tilde{u}, \tilde{v})$ be the two 
   charts of $\OOd$, related by  
$$(\tilde{z}, \tilde{u}, \tilde{v})
   =(1/z,z^d u, z^d v).$$ There is an antiholomorphic involution
   \begin{eqnarray*}
    \tsi\co \OOd&\longrightarrow&\OOd\\
    (z,u,v)&\longrightarrow&(\frac{1}{\bar{z}},
    \bar{z}^d\bar{v},\bar{z}^d\bar{u})
   \end{eqnarray*}
   in terms of the first chart. The involution $\tsi$ covers
   $\si\co \PP\rightarrow\PP$, and the fixed locus of $\tsi$ is $N(d)_\R$.
   Note that $H^0(\PP,\OOd)=0$, and
\[
   H^1(\PP,\mathcal{O}(-d)\oplus\mathcal{O}(-d))=
   \left\{\left.\left(\sum_{j=1}^{d-1}\frac{a_{-j}}{X_0^j X_1^{d-j}},
   \sum_{j=1}^{d-1}\frac{a_j}{X_0^{d-j} X_1^j}\right)
   \;\right|\; a_{-j},a_j\in\C \right\},
\]
   while
\[
   \tsi(a_{-(d-1)},...,a_{-1},a_1,...,a_{d-1})
   =(\bar{a}_{d-1},...,\bar{a}_1,\bar{a}_{-1},...,\bar{a}_{-(d-1)}).
\]
   Therefore, $H^0(D^2,S^1,N(d),N(d)_\R)=0$, and 
\[
   H^1(D^2,S^1,N(d),N(d)_\R)\cong\left\{\left.
   \left(\sum_{j=1}^{d-1}\frac{\bar{a}_j}{z^{d-j}},
   \sum_{j=1}^{d-1}\frac{a_j}{z^j}\right)\;\right|\; a_j\in\C\right\}.
\]
  \end{example}     

 \subsection{Moduli spaces of bordered Riemann surfaces}
 \label{surfacemoduli} Let $M_{g;h}$ denote the moduli space of
 isomorphism classes of bordered Riemann surfaces of type $(g;h)$,
 where $2g+h>2$. It is a double cover of $M(2g+h-1,h,0)$, the moduli
 space of symmetric Riemann surface of type $(2g+h-1, h, 0)$, which is
 a semialgebraic space of real dimension $-3\chi(\Si_{g;h})= 6g+3h-6$
 Sepp{\"a}l{\"a}--Silhol \cite{SS}, where $\chi(\Si_{g;h})=2-2g-h$ is
 the Euler characteristic of a bordered Riemann surface of type
 $(g;h)$.  The covering map is given by $[\Si]\mapsto [\Si_\C]$, and
 the preimage of $[\Si_\C]$ is $\{[\Si], [\overSi]\}$. Note that $\Si$
 and $\overSi$ are not isomorphic bordered Riemann surfaces for
 generic $\Si$, so this cover is generically $2$ to $1$.
   
  \subsection{Stable bordered Riemann surfaces}
  \label{stabbord}

   \begin{df}
    A \emph{node} on a singular bordered Riemann surface $\Si$ is either a 
    singularity on $\inS$ isomorphic to $(0,0)\in\{xy=0\}$ or a singularity
    on $\bS$ isomorphic to $(0,0)\in\{xy=0\}/A$, where $(x,y)$ are coordinates
    on $\C^2$, $A(x,y)=(\bar{x},\bar{y})$ is the complex conjugation. 
    A \emph{nodal bordered Riemann surface} is a singular bordered Riemann
    surface whose singularities are nodes.
   \end{df}   
    
   The notion of morphisms and complex doubles can be easily extended to
   nodal bordered Riemann surfaces. The complex double of a nodal bordered 
   Riemann surface is a nodal compact symmetric Riemann surface.
   The boundary of a nodal Riemann surface is a union of circles, where the
   intersection of any two distinct circles is a finite set.
   
   \begin{df} \label{norm}
    Let $\Si$ be a nodal bordered Riemann surface. The antiholomorphic
    involution $\si$ on its complex double $\Si_\C$ can be lifted to 
    $\hsi\co \widehat{\Si_\C}\rightarrow\widehat{\Si_\C}$, where
    $\widehat{\Si_\C}$ is the normalization of $\Si_\C$ (viewed as a
     complex algebraic curve).
    The {\em normalization} of $\Si=\Si_\C/\bra\si\ket$ is defined to be
    $\hSi=\widehat{\Si_\C}/\bra\hsi\ket$.
   \end{df}
 
    From the above definition, the complex double of the normalization 
    is the normalization of the complex double, ie, 
    $\hSi_\C=\widehat{\Si_\C}$.

   \begin{df}
    A \emph{prestable} bordered Riemann surface is either a smooth bordered
    Riemann surface or a nodal bordered Riemann surface. A \emph{stable}
    bordered Riemann surface is a prestable bordered Riemann surface whose
    automorphism group is finite. 
   \end{df}
 
    The complex double of a stable bordered Riemann surface is a stable
    complex algebraic curve. A smooth bordered Riemann surface of type
    $(g;h)$ is stable if and only if $2g+h>2$.
   
   The stable compactification $\overM_{g;h}$ of $M_{g;h}$ is a double cover
   of the stable compactification $\overM(2g+h-1,h,0)$ of $M(2g+h-1,h,0)$,
   which is compact and Hausdorff \cite{S}. The covering map described in
   \fullref{surfacemoduli} has an obvious extension to 
   $\overM_{g;h}\rightarrow\overM(2g+h-1,h,0)$.
 
\subsection{Riemann--Roch theorem for prestable bordered Riemann surfaces}
  
  Let $\Si$ be a nodal bordered Riemann surface such that
  $H_2(\Si;\Z)=0$.  Let $E$ be a complex vector bundle over $\Si$, and
  let $E_\R$ be a totally real subbundle of $E|_{\bS}$.  In this
  situation, $E$ is a topologically trivial complex
  vector bundle. One can trivialize $E$ and define the
  \emph{Maslov index} $\mu(E,E_\R)$ which is independent of the choice of
  the trivialization as before, since the proof of
  \fullref{Maslov} is valid for nodal $\Si$ as well.

  \begin{df}
   A prestable bordered Riemann surface is \emph{irreducible} if its complex
   double is an irreducible complex algebraic curve.
  \end{df}

  Let $\Si$ be a nodal bordered Riemann surface, let $C_1, ..., C_\nu$
  be the irreducible components of $\Si$ which are (possibly nodal)
  Riemann surfaces, and let $\Si_1, ..., \Si_{\nu'}$ be the remaining 
  irreducible components of $\Si$, which are (possibly nodal) bordered
  Riemann surfaces. 
  Then the irreducible components of $\Si_\C$ are
  \begin{equation}
  \label{components}
   C_1,\ldots,, C_\nu,\ \overline{C}_1,\ldots,\overline{C}_\nu,\ 
   (\Si_1)_\C,\ldots,(\Si_{\nu'})_\C.
  \end{equation}
  Let $E$ be a complex vector bundle over $\Si$
  and let $E_\R$ be a totally real subbundle of $E|_{\bS}$. Observe that
  $H_2(\Si_{i'};\Z)=0$, so $\mu(E|_{\Si_{i'}},E_\R|_{\bS_{i'}})$ is defined for
  $i'=1, ...,\nu'$.

  \begin{df}
   Let $\Si$, $E$, $E_\R$ be as above. The \emph{Maslov index} of $(E,E_R)$ 
   is defined by
\[ 
   \mu(E,E_\R)=2\sum_{i=1}^\nu\deg(E|_{C_i})+
   \sum_{i'=1}^{\nu'}\mu(E|_{\Si_{i'}}, E_\R|_{\bS_{i'}}).
\]
   \end{df}

  Let $\Si$, $C_i$, $\Si_{i'}$ be as above, $i=1,\ldots,\nu$, 
  $i'=1,\ldots,\nu'$. Then
\[ 
  H_2(\Si_\C,\Z)=\bigoplus_{i=1}^\nu \Z[C_i]\oplus 
               \bigoplus_{i=1}^\nu \Z[\overline{C}_i]\oplus
               \bigoplus_{i'=1}^{\nu'} \Z[(\Si_i)_\C]. 
\]
  Set
\[
  [\Si_\C]=\sum_{i=1}^\nu [C_i]+\sum_{i=1}^\nu [\overline{C}_i]+
      \sum_{i'=1}^{\nu'}[(\Si_i)_\C]\in H_2(\Si_\C,\Z).
\]
  With this notation, we have the following generalization of 
  \fullref{degree}:
  \begin{thm}\label{nodaldeg}
   Let $E$ be a complex vector bundle over a prestable bordered Riemann
   surface $\Si$, and let $E_\R$ be a totally real subbundle of $E|_{\bS}$.
   Then there is a complex vector bundle $E_\C$ on $\Si_\C$ together with a 
   conjugate linear involution $\tsi\co E_\C\rightarrow E_\C$ covering the 
   antiholomorphic involution $\si\co \Si_\C\rightarrow\Si_\C$
   such that $E_\C|_\Si=E$ and the fixed locus of $\tsi$
   is $E_\R\rightarrow \bS$. Moreover, we have
\[
   \deg E_\C=\mu(E,E_\R).
\]
  \end{thm}

\proof The construction of $E_\C$ is elementary, 
   as in the smooth case. The reasoning in the proof of 
   \fullref{degree} shows that
\[
\deg(E|_{\Si_{i'}})_\C=\mu(E|_{\Si_{i'}}, E_\R|_{\bS_{i'}}),
\]
so that
$$\eqalignbot{
    \deg E_\C &= \sum_{i=1}^\nu \deg (E_\C|_{C_i})+
                  \sum_{i=1}^\nu \deg (E_\C|_{\Bar{C}_{i}})+
                  \sum_{i'=1}^{\nu'} \deg (E_\C|_{\Si_{i'}})\cr
     &= 2\sum_{i=1}^\nu \deg (E|_{C_i})+\sum_{i'=1}^{\nu'}
        \deg(E|_{\Si_{i'}})_\C\cr
     &= 2\sum_{i=1}^\nu \deg (E|_{C_i})+
         \sum_{i'=1}^{\nu'}\mu(E|_{\Si_{i'}},E_\R|_{\bS_{i'}})\cr
     &= \mu(E,E_\R).}\eqno{\qed}$$
  \begin{df}
   We call the bundle $E_\C$ as in \fullref{nodaldeg} the \emph{topological
   complex double} of $(E,E_\R)$.
  \end{df}
 
  \begin{df}
   Let $E$ be a complex vector bundle of rank $n$ over a prestable bordered
   Riemann surface $\Si$, and let $E_\R$ be a totally real subbundle of 
   $E|_{\bS}$. $(E,E_\R)\rightarrow (\Si,\bS)$ is called a 
   \emph{Riemann--Hilbert bundle} over $\Si$ if there is an open cover
   $\{U_i\mid i\in I\}$ of $\Si$ together with trivializations 
   $\Phi_i\co E|_{U_i}\cong U_i\times \C^n$ such that
   \begin{eqnarray*} 
   \Phi_{ij}\equiv \Phi_i\circ \Phi_j^{-1}\co  
    (U_i\cap U_j)\times \C^n&\longrightarrow&(U_i\cap U_j)\times \C^n\\
    (x,u)&\mapsto&(x, g_{ij}(x)u)
   \end{eqnarray*}
   where $g_{ij}\co U_i\cap U_j\rightarrow GL(n,\C)$ is holomorphic, and
   $g_{ij}(U_i\cap U_j\cap \bS)\subset GL(n,\R)$.  The trivialization
   $\{(U_i,\Phi_i)\mid i\in I\}$ is called a \emph{Riemann--Hilbert
   trivialization} of $(E,E_\R)\rightarrow (\Si,\bS)$.  
\end{df}

   The construction of the \emph{holomorphic complex double} of a 
   Riemann--Hilbert bundle can be easily extended to Riemann--Hilbert
   bundles over nodal bordered Riemann surfaces. 
   Let $(E,E_\R)\rightarrow(\Si,\bS)$ be a Riemann--Hilbert bundle 
   of rank $n$ over a prestable bordered Riemann surface $\Si$, and
   let $E_\C\rightarrow \Si_\C$ be its complex double. 
   We define $Ind(E,E_\R)$ and $\chi(E,E_\R)$ as in
   \fullref{index}. Letting $\tilde{g}$ denote the arithmetic
   genus of $\Si_\C$, we have
   \begin{eqnarray*}
    \chi(E,E_\R)&=&\chi(E_\C)\\
               &=&\deg(E_\C)+n(1-\tilde{g})\\
               &=&\mu(E,E_\R)+n(1-\tilde{g}),
   \end{eqnarray*}
   where the second identity follows from Riemann--Roch and the third
   comes from \fullref{nodaldeg}. Therefore, we have the 
   following Riemann--Roch theorem for prestable bordered Riemann
   surfaces:
 
 \begin{thm}
 \label{nodalRR}
   Let $E$ be a Riemann--Hilbert bundle of rank $n$ over a prestable
   bordered Riemann surface $\Si$. Then
\[
   \chi(E,E_\R)=\mu(E,E_\R)+n(1-\tilde{g}),
\]
   where $\tilde{g}$ is the arithmetic genus of $\Si_\C$, the complex 
   double of $\Si$.
\end{thm}

\section{Stable maps and their moduli spaces}
\label{stable}
  
\subsection{Stable maps} 
\label{mapmoduli}
 
\begin{df} 
\label{stablemap}
Let $\Si$ be a prestable bordered Riemann surface,
let $(X,J,\omega)$ be a symplectic manifold together with an almost complex
structure $J$ compatible with the symplectic form $\omega$, and let
$L\subset X$ be a Lagrangian submanifold of $X$ with respect to $\omega$.
A \emph{prestable map} $f\co (\Si,\bS)\rightarrow (X,L)$ is a continuous 
map which is $J$--holomorphic on $\inS$.
\end{df}
     
\begin{df} \label{mapmor}
\begin{sloppypar}
    A \emph{morphism} between two prestable maps
    $f\co (\Si,\bS)\rightarrow (X,L)$ and $f'\co (\Si',\bS')\rightarrow
    (X,L)$ is a morphism $g\co (\Si,\bS)\rightarrow(\Si',\bS')$ such that
    $f=f'\circ g$.  If $f=f'$ and $g$ is an isomorphism, then $g$ is
    called an automorphism of $f$.  A \emph{stable map\/} is a
    prestable map whose automorphism group is finite.
\end{sloppypar}
\end{df}
  
Let $(X,J,\omega)$ and $L\subset X$ be as in \fullref{stablemap}.
Let 
\[
\MX
\]
denote the moduli space of isomorphism classes of stable maps 
$f\co (\Si,\bS)\rightarrow (X,L)$ such that $f_*[\Si]=\beta\in H_2(X,L;\Z)$ and
$f_*[R_i]=\gamma_i\in H_1(L)$, where $\Si$ is a prestable bordered Riemann 
surface of type $(g;h)$, and $R_1, ..., R_h$ are the connected components 
of $\bS$ with the standard orientations on the $R_i$. Here an
isomorphism between stable maps is an isomorphism in the sense of Definition
\ref{mapmor} which preserves the ordering of the boundary components. 
A necessary condition for the moduli space to be non-empty is
$\partial\beta=\sum\gamma_i$.

The moduli space $\MX$ has extra structure, which will be described elsewhere.
For our present purposes, we content ourselves with remarking that 
$\MX$ should be thought of as an orbi-space due to the presence of nontrivial
automorphisms.  The only additional structure we will use arises from 
deformations and obstructions, to be described in \fullref{slag}.
We will frequently be in the situation where $X$ is a complex algebraic 
manifold, $J$ is the complex structure, and $\omega$ is a K\"ahler form.
In this situation, we will frequently use algebro-geometric language and
related constructions in describing $\MX$.
   
\subsection{Virtual dimension of the moduli space}
\label{slag}

Let $(X,J,\omega)$ be as in \fullref{stablemap}, and let
$L\subset X$ be Lagrangian with respect to $\omega$.  We want to
define and compute the virtual dimension of $\MX$, where $\beta\in
H_2(X,L;\Z)$ and $\gamma_1, ..., \gamma_h\in H_1(L)$. We propose that
there is a tangent-obstruction exact sequence of sheaves on $\MX$ much
as in ordinary Gromov--Witten theory\footnote{The term
$H^0(\Si,\bS,T_\Si, T_{\bS})$ is correct only for smooth $\Si$.  In
general, we need to consider instead vector fields which vanish at the
nodes, together with a contribution for the smoothing of the nodes
\cite{kont}.}
\begin{equation}
\label{toseq}
\begin{array}{ccc}
     0&\rightarrow& H^0(\Si,\bS,T_\Si, T_{\bS})
      \rightarrow H^0(\Si,\bS,f^* T_X,(f|_{\bS})^* T_L) 
      \rightarrow \T^1\\
      &\rightarrow& H^1(\Si,\bS, T_\Si, T_{\bS})
      \rightarrow H^1(\Si, \bS,f^* T_X,(f|_{\bS})^* T_L)
      \rightarrow \T^2\rightarrow 0
\end{array}
\end{equation}
The 4 terms other than the $\T^i$ are labeled by their fibers. 
We will continue to use this description, analogous to that
described in \cite[Section~7.1]{CK} in Gromov--Witten
theory, leaving a more formal treatment for future work.

The virtual (real) dimension of $\MX$ at $f$ is defined to be 
\begin{equation}
\label{vddf}
\mathrm{rank}\left(\T^1\right)-\mathrm{rank}\left(\T^2\right).
\end{equation}
The $\T^i$ need not be bundles, but the rank can be computed fiberwise.
Using \ref{toseq}, the virtual dimension is just
\begin{equation}
\label{vd}
\chi\left(f^* T_X, (f|_{\bS})^* T_L\right)-\chi(T_\Si,T_{\bS}).
\end{equation}
For ease of exposition, assume that $(\Si,\bS)$ is smooth.  By 
\fullref{RR}, we get for the virtual dimension
\[
\mu\left(f^* T_X, (f|_{\bS})^* T_L\right)-\mu\left(\chi(T_\Si,T_{\bS})\right)
+(n-1)\chi(\Sigma),
\]
where $n=\dim X$.
Noting that 
\[
\mu\left(\chi(T_\Si,T_{\bS})\right)=2\chi(\Si)
\]
(as can be seen for example by doubling), the virtual dimension \ref{vd}
becomes 
\begin{equation}
\label{muvd}
\mu\left(f^* T_X, (f|_{\bS})^* T_L\right)+(n-3)\chi(\Si).
\end{equation}
The above discussion leading up to~\ref{muvd} can be extended to general
$\Si$ using \fullref{nodalRR} in place of \fullref{RR}.

Suppose that $L$ is the fixed locus of an antiholomorphic involution
$A\co X\to X$. Then in this situation, $f_\C^* T_X$ is the complex double
of $(f^* T_X,(f|_{\bS})^* T_L)$ and $T_{\Si_{\C}}$ is the complex
double of $(T_\Si,T_{\bS})$ (we only need the complex double of a
bundle, but in this case it happens that the double comes from the
complex double of a map).  This implies that
\[
\mu\left(f^* T_X, (f|_{\bS})^* T_L\right)=\deg\left(f_\C^*T_X\right).
\]
This gives for the virtual dimension
\begin{equation}
\label{vddouble}
\deg\left(f_\C^* T_X\right)+(n-3)\chi(\Si).
\end{equation}
We could have obtained the same result by Riemann--Roch on $\Sigma_\C$,
noting that the genus $\tilde{g}=2g+h-1$ of $\Si_\C$ satisfies
$1-\tilde{g}=\chi(\Sigma)$.  It is clear from \ref{vddouble} that
in this situation, the virtual dimension \ref{vddf} is independent of
the chosen map $f$ in the moduli space.

The formula~\ref{vddouble} coincides with the well known formula 
for the virtual dimension of $\overM_{\tilde{g},0}(X,(f_\C)_*
[\Si_\C])$ in ordinary Gromov--Witten theory \cite[Section 7.1.4]{CK}.
Each of the terms in~\ref{toseq} is in fact the fixed locus of the
induced antiholomorphic involution on the terms of the corresponding
sequence in Gromov--Witten theory.

\begin{rem}
The stack $\overM_{\tilde{g},0}(X,(f_\C)_*[\Si_\C])$ has a natural
antiholomorphic involution induced from that of $X$, and it is natural
to try to compare the fixed locus of this involution with $\MX$.  To
address this, we could have expressed our theory in terms of {\em
symmetric stable maps\/}, ie, stable maps from symmetric Riemann
surfaces to $X$ which are compatible with the symmetry of the Riemann
surface and the antiholomorphic involution on $X$ in the natural
sense.  But there are subtleties: non-isomorphic symmetric stable maps
can give rise to isomorphic stable maps, and \fullref{notdouble}
says in this context that not all symmetric stable
maps arise as doubles of stable maps of bordered Riemann surfaces.
\end{rem}

\begin{rem}
Note that there are subtleties in the notion of the virtual dimension
in the general case, since the Maslov index is a homotopy invariant
rather than a homology invariant in general.
\end{rem}

We will be considering two situations in this paper.  For the first,
consider $(\PP,J,\omega)$, where $J$ is the standard complex
structure, and $\omega$ is the K\"{a}hler form of the Fubini--Study
metric.  Set $\beta=[D^2]\in H_2(\PP, S^1;\Z)$, and $\gamma=[S^1]\in
H_1(S^1)$. We denote $\overM_{g;h}(\PP,S^1|d\beta;n_1\gamma, ...,
n_h\gamma)$ by $\MD$.  

If $f\in\MD$, then $\deg(f_\C^* T_X)=2d$.  So
\ref{vddouble} gives $2(d+2g+h-2)$ as the virtual dimension of
$\MD$.
   
The other situation is $(X,\omega, \Omega)$, where $X$ is a Calabi--Yau
$n$--fold, $\omega$ is the K\"{a}hler form of a Ricci flat metric on
$X$, $\Omega$ is a holomorphic $n$--form, and
\[\frac{\omega^n}{n!}=(-1)^{\frac{n(n-1)}{2}} 
      \left(\frac{i}{2}\right)^n \Omega\wedge \bar{\Omega},
\]
where $\Omega$ is normalized so that it is a calibration
Harvey--Lawson~\cite[Definition 4.1]{HL}. Let $L\subset X$ be a
submanifold which is the fixed locus of an antiholomorphic isometric
involution $A\co X\rightarrow X$ so that $L$ is a special Lagrangian
submanifold of $(X,\omega, \Omega)$.

In this situation, $\deg(f_\C^* T_X)=0$ because $c_1(T_X)=0$. So the virtual 
real dimension of $\MX$ is 
\[
(n-3)\chi(\Sigma)=(n-3)(2-2g-h).
\]
Observe that the virtual dimension does not depend on the classes
$\beta, \gamma_1, ..., \gamma_h$. In particular, the virtual dimension
of $\MX$ is $0$ if $n=3$. 

\begin{rem}
This computation can be done without $A$ and the doubling construction using
the vanishing of the Maslov index for special Lagrangians.
\end{rem}

Continuing to assume the existence of an antiholomorphic involution, we
propose the existence of a virtual fundamental class $[\MX]^\mathrm{vir}$
on $\MX$ of dimension equal to the virtual dimension~\ref{muvd} or 
\ref{vddouble}.  The only assumption we will need to make on this class
is that it can be computed using~\ref{toseq} and a torus action in
our situation, as will be explained in more detail later.

More precisely, we are not proposing the existence of an intrinsic
notion of the virtual fundamental class, but rather a family of
virtual fundamental classes depending on additional choices made at
the boundary of $\MX$.  We will see later that different choices lead
to different computations of the desired enumerative invariants in examples.

\section{Torus action}
\label{ta}

In this section, we define $U(1)$ actions on moduli spaces and compute
the weights of certain $U(1)$ representations that we will need later.

\subsection{Torus action on moduli spaces} 
\label{torus}
 
Let $z$ be an affine coordinate on $\PP$, and put $D^2=\{z: |z|\leq
1\}\subset\PP$, with its boundary circle denoted by $S^1$ with the
standard orientation induced by the complex structure.

Consider the action of $U(1)$ on $(\PP,S^1)$: 
\begin{equation}
\label{action}
z\mapsto e^{-i\theta}z,\qquad e^{i\theta}\in U(1).
\end{equation}
This action has been chosen for consistency with~\cite{CK} and
\cite{FP}.  The action \ref{action} also preserves $D^2$ and
therefore induces an action of $U(1)$ on $\MD$ by composing a stable
map with the action on the image.

More generally, if $(X,L)$ is a pair with a $U(1)$ action, then $U(1)$ 
acts naturally on any space of stable maps to $(X,L)$.

The fixed locus $F_{0;1|d;d}$ of the $U(1)$ action on
$\overM_{0;1}(D^2,S^1|d;d)$ 
consists of a single point corresponding the map 
$f\co D^2\rightarrow D^2$, $x\mapsto x^d$, which has an automorphism group of
order $d$. The fixed locus $F_{0;2|d;n_1,n_2}$ of the $U(1)$ action on 
$\overM_{0;2}(D^2,S^1|d;n_1,n_2)$ consists of a single point corresponding
the map $f\co D_1\cup D_2\rightarrow D^2$, where $D_1$ and $D_2$ are two discs 
glued at the origin to form a node. The restrictions of $f$ to $D_1$ and
to $D_2$ are $x_1\mapsto x_1^{n_1}$ and $x_2\mapsto x_2^{n_2}$, where $x_1$
and $x_2$ are coordinates on $D_1$ and $D_2$, respectively.
For $(g;h)\neq(0;1),(0;2)$, 
the $U(1)$ action on $\MD$ has a single fixed component 
$F_{g;h|n_1,...,n_h}$, consisting of the following stable maps
$f\co (\Sigma,\partial\Sigma)\to (D,S^1)$:
\begin{itemize}
\item The source 
curve $\Sigma$ is a union $\Sigma_0\cup D_1\ldots\cup D_h$, where
$\Sigma_0=(\Sigma_0,p_1,\ldots,p_h)$
is an $h$--pointed stable curve of
genus $g$, and $D_1,\ldots,
D_h$ are discs.  The origin of $D_i$ is glued to $\Sigma_0$ at 
$p_i$ to form a node.
\item $f$ collapses $\Sigma_0$ to 0.
\item $f|_{D_i}\co D_i\to D$ is the cover $f(x)=x^{n_i}$.
\end{itemize}

The fixed locus $F_{g;h|n_1,...,n_h}$ in $\MD$ is an algebraic stack
contained in $\MD$.  It is isomorphic to 
a quotient of $\overM_{g;h}$ by two types of automorphisms:
the automorphisms induced by the covering transformations of the
$f|_{D_i}$, and the automorphisms induced by permutations of the $p_i$
which leave the $n_i$ unchanged.

\subsection{Torus action on holomorphic vector bundles}
\label{torusbund}
 
Let $E_\C$ be a holomorphic vector bundle on $\PP$, and
$\tsi\co  E_\C\rightarrow E_\C$ an antiholomorphic involution which covers 
$\si\co \PP\rightarrow\PP$, given by 
$\si(z)=1/\bar{z}$ in an affine coordinate. The fixed locus of 
$\tsi$ is the total space of a real vector bundle $E_\R\rightarrow S^1$, and
$E_\C|_{S^1}=E_\R\otimes_\R \C$.  The group
$\Cstar$ acts on $\PP$ by $z\rightarrow \lambda z$ for $\lambda\in\Cstar$.
We want to lift the $\Cstar$ action on $\PP$ to $E$ in such a way
that $U(1)\subset\Cstar$
acts on $E_\R$. 

\begin{df}
A lifting of the $\Cstar$ action on $\PP$ to $E$ is
{\em compatible with $\tsi$\/} if
$\tsi(\lambda\cdot v)=
\bar{\lambda}^{-1}\cdot\tsi$ for all $\lambda\in\Cstar$ and $v\in E$.
\end{df}

If the $\Cstar$ action is compatible with $\tsi$, then $U(1)$ automatically
acts on $E_\R$ as desired.
   
We now establish some terminology.  We say that a lifting of a $\Cstar$
action from $\PP$ to a line bundle $L$ has weight $[a,b]$ if the induced
$\Cstar$ actions on the fibers over the fixed points on $\PP$ have
respective weights $a$ and $b$.  We say that a lifting to a sum
$L_1\oplus L_2$ of line bundles has weight $[a,b],\ [c,d]$ if the
lifting to $L_1$ has weight $[a,b]$ and the lifting to $L_2$ has weight
$[c,d]$.

\begin{example} \label{torusbund1}
Let $f\co \PP\rightarrow\PP$ be given in affine coordinates by 
$z=f(x)=x^d$. The antiholomorphic involution 
$\si$ on $\PP$ has a canonical lifting to $T_{\PP}$, thus a canonical 
lifting $\tsi$ to $f^*T_{\PP}$.
The canonical lifting of the torus action on $\PP$ to $f^* T_{\PP}$ has
weights $[1,-1]$ at the respective fixed points \cite[Section 2]{GP} 
and is compatible with $\tsi$. 

The corresponding weights of the $\Cstar$ action on 
$H^0(\PP,f^*T_{\PP})$ are
\[ 
\left\{\frac{d}{d},\frac{d-1}{d}, \cdots,\frac{1}{d}, 0, 
      -\frac{1}{d},\cdots,-\frac{d}{d}\right\}
\] 
with respect to the ordered basis 
\[ 
\mathcal{B}=\left( x^{2d}\frac{\partial}{\partial z}, 
x^{2d-1}\frac{\partial}{\partial z},\cdots, \frac{\partial}{\partial z}\right) 
\]
expressed in affine coordinates. 

We compute that
\begin{equation}
\label{realfields}
\begin{array}{ccl}
H^0(\PP, f^*T_\PP)^{\tsi}&\simeq&
H^0(D^2,S^1,(f|_{D^2})^* T_{\PP},(f|_{S^1})^* T_{S^1})\\
&=&H^0(D^2,S^1, L(2d),L(2d)_\R)
\end{array}
\end{equation}
(see \fullref{double1} for the
definition of $(L(m),L(m)_\R)$).  For later use, note that
\ref{realfields} can be identified with the kernel $V$ of 
\begin{eqnarray*}
 H^0(\PP, \Oh(d))\oplus H^0(D^2,S^1, \C,\R)&\rightarrow&\C\\
 (\sum_{j=0}^d a_j X_0^{d-i} X_1^i,b)&\mapsto& a_d-b.
\end{eqnarray*}
Explicitly, the identification is given by
\begin{equation}
\label{ident}
\begin{array}{ccl}
 H^0(D^2,S^1, L(2d), L(2d)_\R)&\rightarrow&V\\
 \sum_{j=0}^{d-1}(a_j x^j-\bar{a}_j x^{2d-j})+ib x^d&\mapsto&
 (\sum_{j=0}^{d-1} a_j X_0^{d-j}X_1^j+bX_1^d,b) 
\end{array}
\end{equation}
where $a_j\in\C$, $b\in\R$.  

We have a $U(1)$ action on $V$ defined by giving $X_0$ weight
$1/d$ and $X_1$ weight~0.  Then~\ref{ident} preserves the $U(1)$
actions.

We now introduce notation for real representations of $U(1)$. 
We let $(0)_\R$ denote the trivial representation on $\R^1$, and $(w)$ 
denote the representation 
\[
e^{i\theta}\mapsto
\left( \begin{array}{cc}\cos w\theta&-\sin w\theta\\ 
        \sin w\theta&\cos w\theta\end{array}\right)
\]
on $\R^2$ (or $e^{i\theta}\mapsto e^{i w\theta}$ 
on $\C$ if we identify $\R^2$ with $\C$).

With this notation, the corresponding representation of
$U(1)\subset \C^*$ on $V$ is  
\[ 
\left(\frac{d}{d}\right)\oplus
   \left(\frac{d-1}{d}\right)\oplus \cdots
   \left(\frac{1}{d}\right)\oplus(0)_\R.
\]
\end{example}
   
\begin{example}\label{torusbund2}
Let $N=\OO$, and consider the standard antiholomorphic
involution on $N$
\[
\tsi\co N\rightarrow N, \qquad 
   \tsi(z,u,v)=(\frac{1}{\bar{z}}, \bar{z}\bar{v},\bar{z}\bar{u}).
\]
A lifting of
the $\Cstar$ action to $\Oh_{\PP}(-1)$ must have 
weights $[a-1,a]$ at the respective fixed points for some integer
$a$.  

The lifting $([a-1,a],[b-1,b])$ is compatible with $\tsi$ 
if $a+b=1$, ie, if it is of the form $([a-1,a],[-a,1-a])$. 

We now consider such a compatible lifting
with weights $([a-1,a],[-a,1-a])$.
Let $f\co \PP\rightarrow\PP$ be given by $f(x)=x^d$
as in \fullref{torusbund1}. 
The antiholomorphic involution 
on $N$ induces an antiholomorphic involution on 
$f^* N\cong \Oh_{\PP}(-d)\oplus\Oh_{\PP}(-d)$. 
The corresponding representation of 
$\Cstar$ on $H^1(\PP,\Oh_{\PP}(-d)\oplus\Oh_{\PP}(-d))$ has weights
\[ 
\left\{a-\frac{d-1}{d},a-\frac{d-2}{d}, \cdots,
       a-\frac{1}{d}, \frac{1}{d}-a,\cdots,\frac{d-1}{d}-a\right\}
\]
with respect to the ordered basis 
\[       
\begin{array}{c}
\left( \frac{1}{X_0^{d-1} X_1},0 \right),
             \left( \frac{1}{X_0^{d-2} X_1^2},0 \right),
      \ldots,\left( \frac{1}{X_0 X_1^{d-1}},0 \right),\\
             \left( 0,\frac{1}{X_0^{d-1} X_1} \right),
      \ldots,\left( 0,\frac{1}{X_0 X_1^{d-1}} \right),
\end{array}
\]   
where \v{C}ech cohomology has been used (see \cite[page~292]{CK}
for the $a=0$ case).  The orientation on 
\[
\begin{array}{ccc}
H^1(\PP,f^*N)^{\tsi} & \simeq &
H^1(D^2,S^1,(f|_{D^2})^* N,(f|_{S^1})^* N_\R)\\
&=&H^1(D^2,S^1, N(d),N(d)_\R)
\end{array}
\]
is given by identifying it with
the complex vector space $H^1(\PP,\Oh_{\PP}\oplus\Oh_{\PP}(-d))$:
\begin{eqnarray*}
 H^1(D^2,S^1, N(d),N(d)_\R)&\rightarrow&
 H^1(\PP,\Oh_{\PP}\oplus\Oh_{\PP}(-d)) \\
\left(\sum_{j=1}^{d-1}\frac{\bar{a}_j}{z^{d-j}},
  \sum_{j=1}^{d-1}\frac{a_j}{z^j}\right)
&\mapsto& \left(0,\sum_{j=1}^{d-1}\frac{a_j}{X_0^{d-j} X_1^j}\right).
\end{eqnarray*}
Here, $(N(d),N(d)_\R)$ is as in \fullref{double}.2.
The corresponding representation of $U(1)\subset\C^*$ on 
    $H^1(\PP,\Oh_{\PP}\oplus\Oh_{\PP}(-d))$ is 
    \[ \left(\frac{1}{d}-a\right)\oplus
       \left(\frac{2}{d}-a\right)\oplus \cdots
       \left(\frac{d-1}{d}-a\right). \] 
\end{example}

\section{Orientation and the Euler class}
\label{orientation}

Our computations in \fullref{calculation} will use the Euler class
of oriented bundles on moduli spaces of stable maps.  Rather than
attempt to define orientations of the bundles directly, it suffices
for computation to define orientations on their restrictions to the
fixed loci under the $U(1)$ actions.  We presume that the eventual
careful formulation, left for future work, will coincide with the natural
choices which we make here.

The essential point is to understand orientations on bundles on $BU(1)$,
the classifying space for $U(1)$.  We consider these in turn.

The $U(1)$ representation $(w)$ gives rise to a rank~2 bundle on
$BU(1)$.  We abuse notation by denoting this bundle by $(w)$ as well.
Since the $U(1)$ action preserves the standard orientation of $\R^2$,
the bundle $(w)$ inherits a standard orientation.

With this choice, the Euler class of $(w)$ is $w\lambda\in H^2(BU(1),\Z)$,
where $\lambda$ is a generator of $H^2(BU(1),\Z)$.

Similarly, the trivial representation $(0)_\R$ gives rise to the trivial
$\R$--bundle on $BU(1)$, which inherits the standard orientation of $\R$.
The Euler class is 0 in this case.

In the rest of this paper, we will always use these orientations on
$BU(1)$ bundles.  We will also have occasion to consider bundles $B$
on algebraic stacks $F$ with a trivial $U(1)$ action, arising as fixed
loci of a $U(1)$ action on a larger space.  In this case, we have
$F_{U(1)}=F\times BU(1)$, and correspondingly the bundles $B_{U(1)}$
decompose into sums of bundles obtained by tensoring pullbacks of the
above $BU(1)$ bundles with pullbacks of holomorphic bundles on $F$.
Since holomorphic bundles have canonical orientations, we are able to
orient all of our bundles in the sequel.  These orientations will
be used without further comment.

We now show that our choices are compatible with the orientations
defined for $g=0$ in Fukaya--Oh--Ohta--Ono~\cite{FO3}.
Consider the fixed point component
$F_{g;h|n_1,...,n_h}$, isomorphic to a quotient of $\overM_{g,h}$.
Consider a map $f$ in this component described as in \fullref{torus}.
Let $f_i=f|_{D_i}\co D_i\to D^2$ be the multiple
cover map.
The tangent space to moduli space of maps
at this point can be calculated via
\[
T_{\Sigma,p_1,\ldots,p_h}\overM_{g,h}\oplus_i T_{p_i}\Sigma_0\otimes
T_{p_i}D_i
\oplus H^0(\Sigma_i,f_i^* T_D,f_i^*
T_{\partial D}).
\]
All spaces except the last have canonical complex structures, hence canonical
orientations.  The last space can be oriented as in \cite{FO3}.  By pinching
the disc along a circle centered at the origin, a union of a sphere and a
disc are obtained.  A standard gluing argument identifies 
\[
H^0(\Sigma_i,f_i^* T_{D^2},f_i^*
T_{\partial D^2})\simeq \mathrm{Ker}
\left(H^0(\CP^1,\Oh(n_i))\oplus H^0(\mathrm{Disc},\C,\R)\to \C\right).
\]
The key point is that $H^0(\mathrm{Disc},\C, \R)$ is identified
with $\R$ by the orientation on $S^1$ so has a canonical
orientation.  Everything else is complex, so has a canonical orientation.

Using~\ref{ident}, this orientation coincides with the orientation
induced by the $U(1)$ action.  This orientation is used in the proof
of \fullref{contribution}.

The infinitesimal automorphisms are oriented similarly, so the
tangent space to moduli is the quotient space, which is therefore oriented.

\section{Main results}
\label{mainresults}
   
Let $C$ be a smooth rational curve in a Calabi--Yau threefold
$X$ with normal bundle
$N=N_{C/X}=\OO$, where we identify $C$ with $\PP$.

Suppose $X$ admits an antiholomorphic involution $A$ with the
special Lagrangian $L$ as its fixed locus.  Suppose in addition
that $A$ preserves $C$
and $C\cap L=S^1$.
Then $A$ induces an antiholomorphic involution
$A_*\co N\rightarrow N$ which covers $A|_C\co C\rightarrow C$,
and the fixed locus of $A_*$ is $N_\R=N_{S^1/L}$, the
normal bundle of $S^1$ in $L$, a real subbundle of rank two
in $N_{C/X}|_{S^1}$.  
  
We now 
let $(z, u, v)$ and $(\tilde{z}, \tilde{u}, \tilde{v})$ be the two charts
of $\OO$, related by $(\tilde{z},\tilde{u},\tilde{v})=(\frac{1}{z},zu,zv)$.
Let $\tilde{X}$ denote the total space of $\OO$. We assume that under a
suitable identification $N_{C/X}\cong\OO$, the map $A_*$ takes the form
\begin{eqnarray*} 
A_*\co \tilde{X}&\longrightarrow&\tilde{X}\\
   (z,u,v)&\longmapsto&(\frac{1}{\bar{z}},\bar{z}\bar{v},\bar{z}\bar{u})
\end{eqnarray*} 
in terms of the first chart (see \fullref{double2}). The fixed locus
$\tilde{L}$ of $A_*$ gets identified with the total space of
$N_{C\cap L/L}$, and $\tilde{L}$ is a special Lagrangian
submanifold of the noncompact Calabi--Yau 3--fold $\tilde{X}$. Then
$(\tilde{X},\tilde{L})$ can be thought of as a local model for $(X,L)$ near
$C$.

Set $\beta=f_*[D^2]\in H_2(X,L;\Z)$ and $\gamma=f_*[S^1]\in H_1(L)$.
We consider $\MXL$
with $n_1+\ldots+n_h=d$. The virtual dimension of $\MXL$ is $0$.

Given a stable map $f$ in $\MD$,
$i\circ f$ is an element of $\MXL$, so there is an embedding
\[
j\co \MD\rightarrow \MXL,
\] 
whose image is a connected component of $\MXL$ which we denote by
$\MXL_D$.  We will sometimes refer to maps in $\MD$ or $\MXL_D$ as
type $(g;h)$ multiple covers of $D^2$ of degree $(n_1,\ldots,n_h)$.
We have
\[
\dim[\MD]^\mathrm{vir}=2(d+2g+h-2)
\]
by~\ref{vddouble}.  Let
\[
\pi\co  \U\rightarrow \MD
\]
be the universal family of stable maps, 
$\mu\co \U\rightarrow D^2$ be the evaluation map, 
and let $\N$ be the sheaf of local holomorphic sections of $N|_{D^2}$
with boundary values in $N_\R$. 

Now let $f\co (\Si,\bS)\to(D^2,S^1)$ be a stable map with 
$$f_*[\Sigma]=d[D^2]
\in H^2(D^2,S^1,\Z).$$  
Then 
\begin{eqnarray*}
H^0(\Si,f^*N,f^*N_\R)&=&0\\
\dim H^1(\Si,f^*N,f^*N_\R)&=&2(d+2g-2+h).
\end{eqnarray*}
Thus the  sheaf $R^1 \pi_* \mu^* \N$ is a rank $2(d+2g+h-2)$ 
real vector bundle over $\MD$, which we have called the obstruction bundle.

Note the diagram
\begin{equation}
\label{compobs}
\begin{array}{ccc}
&&0\\
&&\downarrow\\
H^i(T_\Si,T_{\bS})&\to&H^i(f^*T_{D^2},f^*T_{S^1})\\
\downarrow&&\downarrow\\
H^i(T_\Si,T_{\bS})&\to&H^i(f^*T_X,f^*T_L))\\
&&\downarrow\\
&&H^i(f^*N_\C,f^*N_\R)\\
&&\downarrow\\
&&0
\end{array}
\end{equation}
Combining~\ref{toseq} (for $(X,L)$ and for $(D^2,S^1)$) and \ref{compobs},
we are led to the following claim:

The virtual fundamental classes are related by
\begin{equation}
\label{vfformula}
\begin{array}{c}
j_*\left(e(R^1\pi_*\mu^*\N)\cap
[\overM_{g,h}(D^2,S^1\mid d;n_1,\ldots,n_h)]^{\mathrm{vir}}\right)
=\\
\ [\MXL_D]^{\mathrm{vir}}.
\end{array}
\end{equation}
We leave a more precise formulation of this claim and its proof for 
future work.

Therefore,%
\footnote{Recall that we mentioned at the end of \fullref{stable}
that the virtual fundamental class depends on choices.  In the following
formula and the rest of this section, we are actually fixing a particularly
simple choice, corresponding to a specific torus action, in order to
compare with \cite{OV}.  The general
case and its relationship to this section will be discussed in
 \fullref{outline}.}
   \[C(g;h|d;n_1, ...,n_h)=\int_{\left[\MD\right]^{\mathrm{vir}}}e(R^1\pi_*\mu^*\N)\]
   is the contribution to the type $(g;h)$ enumerative invariants
   of $(X,L)$ from multiple covers of $D^2$ of
degree $(n_1, ..., n_h)$.

Alternatively, $C(g;h|d;n_1, ..., n_h)$ can be viewed as enumerative
invariants of $(\tilde{X},\tilde{L})$. 

Since for a K\"ahler class $\omega$ on $X$ we have
\[
\int_{D^2}f^*\omega =\frac12\int_{\PP}f^*\omega
\]
(extending $f$ to $\PP$ as its complex double, see \fullref{cxdm}), 
the truncated type
$(g;h)$ prepotential takes the form
\[
    F_{g;h}(t,y_1, ..., y_h)=
\sum_{n_1+\cdots+n_h=d\atop n_1,...,n_h\geq 0}
 C(g;h|d;n_1, ..., n_h)
    y_1^{n_1}\cdots y_h^{n_h} e^{-\frac{dt}{2}}
\]
where we have put $d=n_1+\ldots+n_h$.
The truncated all-genus potential is 
\begin{eqnarray*}  
F(\lambda,t, y_1,y_2, ...)&=&\sum_{g=0}^\infty\sum_{h=1}^\infty
  \lambda^{-\chi(\Si_{g;h})}F_{g;h}(t,y_1, ..., y_h)\\
  &=&\sum_{g=0}^\infty\sum_{h=1}^\infty
     \lambda^{2g+h-2}F_{g;h}(t,y_1, ..., y_h)
\end{eqnarray*}
Our final assumption is that a natural extension of 
the virtual localization formula of \cite{GP} holds in
this context.\footnote{In particular, this extension requires that
the torus action is chosen to preserve the additional structure, as
will be described in \fullref{outline}.}
  Rather than formulate this assumption in generality,
we will illustrate it in our situation.
The generalization is straightforward.

\begin{propo}
\label{contribution}
Under our assumptions, we have
\begin{eqnarray*}
C(g;1|d;d)&=&d^{2g-2}b_g\\
C(g;h|n_1+\cdots n_h,n_1,...,n_h)&=&0\;\;\;\textup{ for }h>1.
\end{eqnarray*}
\end{propo}

We will prove \fullref{contribution} in
\fullref{calculation}.  The next result follows immediately.

\begin{tmm}
\label{ovresult}
\[
F(\lambda,t,y)=\sum_{d=1}^\infty\frac{e^{-\frac{dt}{2}}y^d}
            {2d\sin\frac{\lambda d}{2}}.
\]
\end{tmm}

This is the multiple cover formula for the disc,
first obtained obtained by string duality in \cite{OV}.

\proof
We compute
\begin{eqnarray*}
F_{g;1}(t,y)&=&\sum_{d=1}^\infty C(g;1|d;d)y^d e^{-\frac{dt}{2}}
        =b_g\sum_{d=1}^\infty d^{2g-2}y^d e^{-\frac{dt}{2}}\\
F_{g;h}(t,y)&=&0\;\;\; \textup{ for }h>1
\end{eqnarray*}
Then
$$\eqalignbot{
    F(\lambda,t,y_1,y_2,...)&=\sum_{g=0}^\infty\lambda^{2g-1}F_{g;1}(t,y)\cr
    &=\sum_{g=0}^\infty \lambda^{2g-1}b_g\sum_{d=1}^\infty 
       d^{2g-2}y^d e^{-\frac{dt}{2}}\cr
    &=\sum_{d=1}^\infty\frac{y^d e^{-\frac{dt}{2}}}{\lambda d^2}
        \sum_{g=0}^\infty b_g(\lambda d)^{2g}\cr
    &=\sum_{d=1}^\infty\frac{y^d e^{-\frac{dt}{2}}}{\lambda d^2}
        \frac{\frac{\lambda d}{2}}{\sin\frac{\lambda d}{2}}\cr
    &=\sum_{d=1}^\infty\frac{e^{-\frac{dt}{2}}y^d}
            {2d\sin\frac{\lambda d}{2}}.}\eqno{\qed}
$$

\section{Final calculations}
\label{calculation}
   
\subsection{Outline}
\label{outline}

In this section, we perform our main calculation and use the result to 
prove \fullref{contribution}.
We want to calculate 
\begin{equation}
\label{theinvt}
\int_{\left[\MD\right]^{\mathrm{vir}}}e(R^1\pi_*\mu^*\N)
\end{equation}
by localization. As before, let $\Cstar$ act on $\OO$ with weights
$([a-1,a], [-a,1-a])$.  From
\fullref{torus} we know that there is only one fixed component.
This lifting is compatible with the antiholomorphic
involution on $\OO$, as explained in \fullref{torusbund2}.  

Recall that over $z=e^{i\theta}\in\partial D^2$, we have $(N_\R)_z=
\{(u,v)\in\C^2|v=\bar{z}\bar{u}\}$.  The additional data needed to
define our invariants turns out to be completely determined by the
topological class of a real 1 dimensional subbundle of $N_\R$, just
like the topological class of a real 1 dimensional subbundle of a
complex line bundle on a bordered Riemann surface is needed to define
its generalized Maslov index.  In order for localization to be valid,
we propose that the torus action is required to preserve this
subbundle.

We define the subbundle $N_a\subset N_\R$ as the subbundle characterized by
\[
(N_a)_z=\{(e^{i\theta(a-1)}r,e^{-i\theta a}r)|r\in\R\}.
\]
Note that the $U(1)$ action on $N$ with weights $([a-1,a],[-a,1-a])$ is
the only compatible lifting of the torus action on $\PP$ which preserves
$N_a$.  Note also that the $\{N_a\}$ represent all of the topological 
classes of rank 1 subbundles of $N_\R$.

We denote by $C(g,h|d;n_1,\ldots,n_h|a)$ 
the invariant \ref{theinvt} to emphasize
its dependence on $a$.  The invariant $C(g,h|d;n_1,\ldots,n_h)$ of
\fullref{mainresults} is the special case $a=0$ of \ref{theinvt}.

\subsection{The obstruction bundle} 
\label{obstruction}

Let $f\co \Si \rightarrow \PP$ be a generic stable map in $\fix$, where
$(g;h)\neq(0;1),(0;2)$. Write $\Si=\Si_0\cup D_1 \cup \cdots \cup D_h$,  where
$\Si_0$ is a curve of genus $g$, and $D_i$ is a disc for $i=1, ..., h$. Let
$z_i$ be the coordinate on $D_i$, identifying $D_i$ with
$D^2\subset \C$. The point $z_i=0$ is identified with $p_i\in \Si_0$,
where $p_1, ..., p_h$ are distinct points on $\Si_0$, forming nodes on the
nodal bordered Riemann surface $\Si$.
Note that $(\Si_0,p_1, ..., p_h)$ represents an element in $\overM_{g,h}$. 
The restriction of $f$ on $D_i$ is $z_i\mapsto z_i^{n_i}$.
Conversely, all $f\co \Sigma\to\PP$ as described above are elements
of $F_{g;h\mid d;n_1,\ldots,n_h}$.

Let $\pi\co  \U\rightarrow \MD$ be the universal family of stable maps, 
let $\mu\co \U\rightarrow D^2$ be the evaluation map, 
and let $\N$ be the sheaf of local holomorphic sections of $N|_{D^2}$
with boundary values in $N_\R$, as in \fullref{mainresults}.
The fiber of the obstruction bundle $R^1\pi_* \mu^* \N$ at 
$[f]$ is $H^1(\Si, f^*\N)=H^1(\Si,\partial \Si,f^* N,f^* N_\R)$.
Consider the normalization sequence
\[ 
0\rightarrow f^*\N
\rightarrow \bigoplus_{i=1}^h (f|_{D_i})^*\N \oplus (f|_{\Si_0})^*\N 
\rightarrow \bigoplus_{i=1}^h (f^*\N)_{p_i}\rightarrow 0.
\]
The corresponding long exact sequence reads
\begin{equation}
\label{normcoh}
\begin{array}{ccc}
0&\rightarrow& H^0(\Si_0,(f|_{\Si_0})^*\N) 
\rightarrow \bigoplus_{i=1}^h H^0((f^*\N)_{p_i})
\rightarrow H^1(\Si,f^*\N)\\
&\rightarrow& \bigoplus_{i=1}^h H^1(D_i,(f|_{D_i})^*\N)
\oplus H^1(\Si_0, (f|_{\Si_0})^*\N)\rightarrow 0 
\end{array}
\end{equation}
 The vector spaces  $H^0(\Si_0,(f|_{\Si_0})^*\N)$,
$H^0((f^*\N)_{p_i})$, and  $H^1(D_i,(f|_{D_i})^*\N)$ appearing in
\ref{normcoh}
fit together to form trivial complex vector bundles over the component
$\fix$ which
are not necessarily trivial as  equivariant vector bundles.
Let $\Cstar$ act on $\OO$ with weights $([a-1,a], [-a,1-a])$.
Then the weights of $U(1)$ in \ref{normcoh} are (see 
\fullref{torusbund2})
   \begin{eqnarray*}
    0&\rightarrow (a-1)\oplus(-a)
       \rightarrow \bigoplus_{i=1}^h ((a-1)\oplus(-a))
       \rightarrow H^1(\Si,f^* \N)\\
      &\rightarrow H^1(\Si_0,\Oh_{\Si_0})\otimes((a-1)\oplus(-a))\oplus\\
&\bigoplus_{i=1}^h
        \left( \left(\frac{1}{n_i}-a\right)\oplus
               \left(\frac{2}{n_i}-a\right)\oplus\cdots\oplus
               \left(\frac{n_i-1}{n_i}-a\right) \right)
\rightarrow 0.
   \end{eqnarray*}
The map
$$i_{g;h|d;n_1,\ldots,n_h}\co \overM_{g,h}\rightarrow  \MD$$ 
has degree 
$n_1\cdots n_h$ onto its image
$F_{g;h|d;n_1,\ldots,n_h}$, and
\[ 
\begin{array}{l}
i_{g;h|d;n_1,\ldots,n_h}^* R^1\pi_*\mu^*\N= \\
\hskip.1truein\left(\hodge^\vee \otimes ((a-1)\oplus(-a))\right)
   \oplus\bigoplus_{i=1}^h\prod_{j=1}^{n_i-1}\left(\frac{j}{n_i}-a \right)
   \oplus\bigoplus_{i=1}^{h-1}((a-1)\oplus(-a))
\end{array}
\] 
where $\hodge$ is the Hodge bundle over $\overM_{g,h}$. Therefore,
\begin{equation}
\label{virtfgen}
\begin{array}{ll}
&e_{U(1)}(i_{g;h|d;n_1,\ldots,n_h}^* R^1\pi_*\mu^*\N)\\
=&c_g(\hodge^\vee((a-1)\lambda))c_g(\hodge^\vee(-a\lambda))
\lambda^{d+h-2}(a(1-a))^{h-1}
   \prod_{i=1}^h\frac{\prod_{j=1}^{n_i-1}(j- n_i a)}{n_i^{n_i-1}}
\end{array}
\end{equation}
where $\lambda_g=c_g(\hodge)$, and $\lambda\in H^2(BU(1),\Z)$ is the
generator discussed in \fullref{orientation}.
  
We now consider the case $(g;h)=(0;1)$. Then $F_{0;1|d;d}$ consists of 
a single point corresponding to the map $f\co D^2\rightarrow D^2$, 
$x\mapsto x^d$. In this case,
\[
H^1(D^2, f^*\N)= \left(\frac{1}{d}-a\right)\oplus
            \left(\frac{2}{d}-a\right)\oplus\cdots\oplus
            \left(\frac{d-1}{d}-a\right).
\]
This gives
\begin{equation}
e_{U(1)}(i_{0,d}^* R^1\pi_*\mu^*\N)=
\lambda^{d-1}\frac{\prod_{j=1}^{d-1}(j-da)}{d^{d-1}}.
\end{equation}
We finally consider the case $(g;h)=(0;2)$. Then $F_{0;2|d;n_1,n_2}$
consists of a single point corresponding to the map $f\co D_1\cup
D_2\rightarrow D^2$.  Here $D_1$ and $D_2$ are identified at the
origin, and the restrictions of $f$ to $D_1$ and to $D_2$ are
$x_1\mapsto x_1^{n_1}$ and $x_2\mapsto x_2^{n_2}$, where $x_1$ and
$x_2$ are coordinates on $D_1$ and $D_2$, respectively.  In this case,
\[
H^1(D^2, f^*\N)= (a-1)\oplus(-a)\bigoplus_{i=1}^2\bigoplus_{j=1}^{n_i-1}
   \left(\frac{j}{d}-a\right).
\]
This gives
\begin{equation}
e_{U(1)}(i_{0,d}^* R^1\pi_*\mu^*\N)=
\lambda^d a(1-a)\prod_{i=1}^2\frac{\prod_{j=1}^{n_i-1}(j-n_i a)}{n_i^{n_i-1}}.
\end{equation}
We conclude that \ref{virtfgen} is valid for $(g;h)=(0;1),(0;2)$ as well.

\subsection{The virtual normal bundle}
\label{virtual}
    
We next compute the equivariant Euler class of the virtual normal
bundle of the fixed point component $F_{g;h|d;n_1,\ldots,n_h}$ in
$\MD$. The computation is very similar to that in \cite[Section
4]{GP}.
 
There is a tangent-obstruction exact sequence of
sheaves on  $F_{g;h|d;n_1,\ldots,n_h}$:
\begin{equation} 
\label{toxg}
\begin{array}{ccccccccc}
 0&\rightarrow& H^0(\Si,\bS,T_\Si, T_{\bS})
   \rightarrow H^0(\Si,\bS,f^* T_{\PP},(f|_{\bS})^* T_{S^1}) 
   \rightarrow \T^1\\
  &\rightarrow& H^1(\Si,\bS, T_\Si, T_{\bS})
   \rightarrow H^0(\Si, \bS,f^* T_{\PP},(f|_{\bS})^* T_{S^1})
   \rightarrow \T^2\rightarrow 0
\end{array}
\end{equation}
The 4 terms in \ref{toxg} other than the sheaves $\T^i$ are vector
bundles and are labeled by fibers. We use \cite{GP} and the notation
there, so that $\T^i=\T^{i,f}\oplus\T^{i,m}$, where $T^{i,f}$ is the
fixed part and $\T^{i,m}$ is the moving part of $\T^i$ under the
$U(1)$ action. Then the $\T^{i,m}$ will determine the virtual normal
bundle of $F_{g;h|d;n_1,\ldots,n_h}$ in $\MD$:
\[ 
e_{U(1)}(N^{\mathrm{vir}})=
\frac{e_{U(1)}(B_2^m)e_{U(1)}(B_4^m)}{e_{U(1)}(B_1^m)e_{U(1)}(B_5^m)},
\]
where the $B_i^m$ denote the moving part of the $i$th term in the
tangent-obstruction exact sequence \ref{toxg}.

We first consider the cases $(g;h)\neq(0;1),(0;2)$. We have 
$B_1=B_1^f=\oplus_{i=1}^h (0)_\R$, which
corresponds to the rotations of the $h$ disc components of the domain curve,
so $e_{U(1)}(B_1^m)=1$. Next, $B_4^m=\oplus_{i=1}^h T_{p_i}\Si_0\otimes T_0
D_i$  corresponds to deformations of the $h$ nodes of the domain. To compute
this, note that $T^*_{p_i}\Si_0$ is the fiber
of $\mathbb{L}_i=s_i^*\omega_\pi\rightarrow \overM_{g,h}$ at 
$[(\Si_0,p_1,\ldots,p_h)]$,
where $\omega_\pi$ is the relative dualizing sheaf of the universal
curve over $\overM_{g,h}$, and $s_i$ is the section corresponding to 
the $i$th marked point. Since $T_0 D_i=(1/n_i)$, we get 
\[
e_{U(1)}(B_4^m)=\prod_{i=1}^h 
\left( c_1(\mathbb{L}_i^\vee)+\frac{\lambda}{n_i}\right)
=\prod_{i=1}^h\left(\frac{\lambda}{n_i}-\psi_i\right),
\]
where $\psi_i=c_1(\mathbb{L}_i)$.
 
For $B_2^m$ and $B_5^m$, consider the normalization exact sequence
\begin{eqnarray*}
0&\rightarrow& (f^* T_{\PP},(f|_{\partial\Si})^* T_{S^1})  
  \rightarrow (f|_{\Si_0})^* T_{\PP}\oplus\bigoplus_{i=1}^h 
    \left( (f|_{D_i})^* T_{\PP},(f|_{\partial D_i})^* T_{S^1}\right)\\
  &\rightarrow& \bigoplus_{i=1}^h(T\PP)_0\rightarrow 0.
\end{eqnarray*}
The corresponding long exact sequence reads
\begin{eqnarray*}
 0&\rightarrow&H^0(\Si,\bS, f^* T_{\PP},(f|_{\partial\Si})^* T_{S^1})\\
   &\rightarrow& H^0(\Si_0,(f|_{\Si_0})^* T_{\PP})\oplus\bigoplus_{i=1}^h 
H^0(D_i,\partial D_i,(f|_{D_i})^* T_{\PP},(f|_{\partial D_i})^* T_{S^1})
\\
 & \rightarrow& \bigoplus_{i=1}^h (T\PP)_0 
\rightarrow H^1(\Si,\bS, f^* T_{\PP},(f|_{\partial\Si})^* T_{S^1})  
  \rightarrow H^1(\Si_0,(f|_{\Si_0})^* T_{\PP})\\
&\rightarrow& 0.
\end{eqnarray*}
The representations of $U(1)$ are (see \fullref{torusbund1})
\begin{equation}
\label{towt}
\begin{array}{ccl}
 0&\rightarrow&H^0(\Si,\bS, f^* T_{\PP},(f|_{\partial\Si})^* T_{S^1})\\
 &\rightarrow& H^0(\Si_0,\Oh_{\Si_0})\otimes(1)\oplus
 \bigoplus_{i=1}^h\left(\bigoplus_{j=1}^{n_i}
 \left(\frac{j}{n_i}\right)\oplus (0)_\R\right)\\
  &\rightarrow&\bigoplus_{i=1}^h (1)
   \rightarrow H^1(\Si,\bS, f^* T_{\PP},(f|_{\partial\Si})^* T_{S^1})  
  \rightarrow H^1(\Si_0,\Oh_{\Si_0})\otimes(1)\rightarrow 0.
\end{array}
\end{equation}
The trivial representation lies in the fixed part, so \ref{towt} gives
\[
\frac{e_{U(1)}(B_2^m)}{e_{U(1)}(B_5^m)}=
    \frac{\lambda^{d+1-h}}{c_g(\hodge^\vee(\lambda))}
    \prod_{i=1}^h\frac{n_i !}{n_i^{n_i}}.
\]
Therefore 
\begin{equation}
\label{virtngen}
e_{U(1)}(N^{\mathrm{vir}})=\frac{\lambda^{d+1-h}}{c_g(\hodge^\vee(\lambda))}
\prod_{i=1}^h\left(
\frac{n_i !}{n_i^{n_i}}\left(\frac{\lambda}{n_i}-\psi_i\right)\right).
\end{equation}
We now consider the case $(g;h)=(0;1)$. The representations of $U(1)$ in
\ref{toxg} are
\[ 
0\rightarrow\left(\frac{1}{d}\right)\oplus(0)_\R
 \rightarrow\bigoplus_{i=1}^{d}\left(\left(\frac{i}{d}\right)\oplus(0)_\R\right)
 \rightarrow\T^1\rightarrow 0\rightarrow 0\rightarrow\T^2\rightarrow 0. 
\] 
Therefore
\begin{equation}
\label{virtn0}
e_{U(1)}(N^{\mathrm{vir}})=\lambda^{d-1}\frac{d!}{d^{d-1}}.
\end{equation}
We finally consider the case $(g;h)=(0;2)$. We have 
$B_1=B_1^f=(0)_\R\oplus(0)_\R$, which corresponds to the rotations
of the two disc components, so $e_{U(1)}(B_1^m)=1$. Next, 
$B_4^m=T_0 D_1\otimes T_0 D_2$ corresponds to deformation of the node
of the domain. Since $T_0 D_i=\left(\frac{1}{n_i}\right)$, we get
$e_{U(1)}(B_4^m)=\left(\frac{1}{n_1}+\frac{1}{n_2}\right)\lambda$.

For $B_2^m$ and $B_5^m$, consider the normalization exact sequence
\[
0\rightarrow (f^* T_{\PP},(f|_{\partial\Si})^* T_{S^1})  
  \rightarrow \bigoplus_{i=1}^2 
    \left( (f|_{D_i})^* T_{\PP},(f|_{\partial D_i})^* T_{S^1}\right)
  \rightarrow (T\PP)_0\rightarrow 0.
\]
The corresponding long exact sequence reads
\begin{eqnarray*}
 0&\rightarrow&H^0(\Si,\bS, f^* T_{\PP},(f|_{\partial\Si})^* T_{S^1})
  \\
&\rightarrow& \bigoplus_{i=1}^2 H^0(D_i,\partial D_i,
   (f|_{D_i})^* T_{\PP},(f|_{\partial D_i})^* T_{S^1})\\
 &\rightarrow& (T\PP)_0
 \rightarrow H^1(\Si,\bS, f^* T_{\PP},(f|_{\partial\Si})^* T_{S^1})  
  \rightarrow 0.
\end{eqnarray*}
The representations of $U(1)$ are (see \fullref{torusbund1})
\begin{equation}\label{towt2}
\begin{array}{ccl}
 0&\rightarrow&H^0(\Si,\bS, f^* T_{\PP},(f|_{\partial\Si})^* T_{S^1})
\rightarrow \bigoplus_{i=1}^2 \left(\bigoplus_{j=1}^{n_i}
 \left(\frac{j}{n_i}\right)\oplus (0)_\R\right)\\
  &\rightarrow&(1)
   \rightarrow H^1(\Si,\bS, f^* T_{\PP},(f|_{\partial\Si})^* T_{S^1})  
  \rightarrow 0.
\end{array}
\end{equation}
The trivial representation lies in the fixed part, so \ref{towt2} gives
\[
\frac{e_{U(1)}(B_2^m)}{e_{U(1)}(B_5^m)}=\lambda^{d-1}
\frac{n_1!n_2!}{n_1^{n_1}n_2^{n_2}}.
\]
Therefore 
\begin{equation}
e_{U(1)}(N^{\mathrm{vir}})=\lambda^d
\left(\frac{1}{n_1}+\frac{1}{n_2}\right)
\frac{n_1!n_2!}{n_1^{n_1}n_2^{n_2}}.
\end{equation}

\subsection{Conclusion}               
From our calculations in Sections \ref{obstruction} and \ref{virtual},
we can now compute the sought-after invariants.
\begin{eqnarray*}
 C(0;1|d;d|a)&=&\frac{1}{d}\int_{\mathrm{pt}_{U(1)}}
 \frac{i_{0;1|d;d}^* e_{U(1)}(R^1\pi_*\mu^*\N)}{e_{U(1)}(N^{\mathrm{vir}})}\\
 &=&\left(\frac{\prod_{j=1}^{d-1}(j-da)}{(d-1)!}\right)
 \frac{1}{d^2},\\
 C(0;2|d;n_1,n_2|a)&=&\frac{1}{n_1 n_2}\int_{\mathrm{pt}_{U(1)}}
 \frac{i_{0;2|d;n_1,n_2}^* e_{U(1)}(R^1\pi_*\mu^*\N)}
{e_{U(1)}(N^{\mathrm{vir}})}\\
 &=&a(1-a)\left(\prod_{i=1}^2\frac{\prod_{j=1}^{n_i-1}(j-n_i a)}{(n_i-1)!}
    \right)\frac{1}{d}.
\end{eqnarray*}
Here $\mathrm{pt}_{U(1)}$ denotes the ``equivariant point'', which is
isomorphic to $BU(1)$.

For $(g;h)\neq(0;1),(0;2)$, we get
\begin{eqnarray*}
&&C(g;h|d;n_1,\ldots,n_h|a)\\
&=&\frac{1}{n_1\cdots n_h}
 \int_{{\overM_{g,h}}_{U(1)}}\frac{i_{g;h\mid d;n_1,\ldots,n_h}^*
      e_{U(1)}(R^1\pi_*\mu^*\N)}{e_{U(1)}(N^{\mathrm{vir}})}\\
&=&(a(1-a))^{h-1}
 \left(\prod_{i=1}^h\frac{\prod_{j=1}^{n_i-1}(j-n_i a)}{(n_i-1)!}\right)\cdot\\
&& \int_{\overM_{g,h}}\frac{c_g(\hodge^\vee(\lambda))
 c_g(\hodge^\vee((a-1)\lambda))c_g(\hodge^\vee(-a\lambda))\lambda^{2h-3}}
 {\prod_{i=1}^h (\lambda-n_i\psi_i)},
\end{eqnarray*}
where ${\overM_{g,h}}_{U(1)}$ denotes $\overM_{g,h}\times BU(1)$.
In particular, we get for $h\geq 3$,
\begin{eqnarray*}
&&C(0;h|d;n_1,\ldots,n_h|a)\\
&=& (a(1-a))^{h-1}
  \left(\prod_{i=1}^h\frac{\prod_{j=1}^{n_i-1}(j-n_i a)}{(n_i-1)!}\right)
    \int_{{\overM_{0,h}}_{U(1)}}\frac{\lambda^{2h-3}}
    {\prod_{i=1}^h (\lambda-n_i\psi_i)}
\end{eqnarray*}

\begin{eqnarray*}
&=& (a(1-a))^{h-1}
  \left(\prod_{i=1}^h\frac{\prod_{j=1}^{n_i-1}(j-n_i a)}{(n_i-1)!}\right)\cdot
\\
&&    \sum_{\begin{array}{c}k_1+\ldots k_h=h-3\\k_1,\ldots,k_h\geq 0
          \end{array}} n_1^{k_1}\cdots n_h^{k_h}
     \int_{\overM_{0,h}}\psi_1^{k_1}\cdots\psi_h^{k_h}.
\end{eqnarray*}
But
\[
\int_{\overM_{0,h}}\psi_1^{k_1}\cdots\psi_h^{k_h}
=\left(\begin{array}{c}h-3\\k_1,\ldots,k_h\end{array} \right),
\]
so 
\begin{eqnarray*}
&&C(0;h|d;n_1,\ldots,n_h|a)\\
&=& (a(1-a))^{h-1}
   \left(\prod_{i=1}^h\frac{\prod_{j=1}^{n_i-1}(j-n_i a)}{(n_i-1)!}\right)\cdot
\\
&&    \sum_{\begin{array}{c}k_1+\ldots k_h=h-3\\k_1,\ldots k_h\geq 0
          \end{array}} n_1^{k_1}\cdots n_h^{k_h}
     \left(\begin{array}{c}h-3\\k_1,\ldots,k_h\end{array} \right)\\
&=& (a(1-a))^{h-1}
   \left(\prod_{i=1}^h\frac{\prod_{j=1}^{n_i-1}(j-n_i a)}{(n_i-1)!}\right)
    (n_1+\ldots+n_h)^{h-3}\\
&=& (a(1-a))^{h-1}
   \left(\prod_{i=1}^h\frac{\prod_{j=1}^{n_i-1}(j-n_i a)}{(n_i-1)!}\right)
    d^{h-3}.
\end{eqnarray*}
Observe the symmetry
\[ 
C(g;h|d;n_1,\ldots,n_h|1-a)=(-1)^{d-h}C(g;h|d;n_1,\ldots,n_h|a),
\]
which is a consequence of the symmetry of $\Oh_{\PP}(-1)\oplus\Oh_{\PP}(-1)$
obtained from interchanging the two factors, together with our orientation
choices (or it can be checked directly).
It is also of interest to observe that for $a>0$, we have the identity
\[
\frac{\prod_{j=1}^{n_i-1}(j-n_i a)}{(n_i-1)!}=(-1)^{n_i-1}
\left(\begin{array}{c}n_i a -1\\n_i -1\end{array}\right)
\]
which is an integer.  There is a similar formula for $a\le0.$

We may summarize our results as follows:

\begin{pro}
Let $a$ be a positive integer. Then
\begin{eqnarray*}
&&(-1)^{d-h}C(0;h|d;n_1,\ldots,n_h|a)=C(0;h|d;n_1,\ldots,n_h|1-a)\\
&=&(a(1-a))^{h-1}
   \prod_{i=1}^h\left(\begin{array}{c}n_i a -1\\n_i-1 \end{array}\right)d^{h-3}.
\end{eqnarray*}
For $g>0$,
\begin{eqnarray*}
&&(-1)^{d-h}C(g;h|d;n_1,\ldots,n_h|a)=C(g;h|d;n_1,\ldots,n_h|1-a)\\
&=&(a(1-a))^{h-1}
 \prod_{i=1}^h\left(\begin{array}{c}n_i a -1\\n_i-1\end{array}\right)\cdot\\
&& \int_{{\overM_{g,h}}_{U(1)}}\frac{c_g(\hodge^\vee(\lambda))
 c_g(\hodge^\vee((a-1)\lambda))c_g(\hodge^\vee(-a\lambda))\lambda^{2h-3}}
 {\prod_{i=1}^h (\lambda-n_i\psi_i)}.
\end{eqnarray*}
\end{pro}

Our results for $g=0$ agree with results of \cite{AKV} up to a sign (eg,
up to a choice of orientation).  Our conventions differ from theirs by a 
sign of $(-1)^{ad}$ (our ``$a$'' is their ``$p$'').  They also note integer
invariants $N_d$ depending on $a$, yielding the identity
\begin{equation}
\label{nd}
(-1)^{ad}C(0,1|d,d|a)=\sum_{k|d}\frac{N_{d/k}}{k^2}.
\end{equation}
In the case $a=0$, we just have $N_1=1$ and $N_d=0$ for all $d\ge1$, in 
which case \ref{nd} is just the multiple cover formula for the disc.  In
the general case, \ref{nd} can be viewed as multiple cover formula with
different integer invariants, in the same sense as the formula of \cite{gov}
for stable maps of closed Riemann surfaces.

\proof[Proof of \fullref{contribution}]
Setting 
\[
C(g;h|d;n_1,...,n_h)=C(g;h|d;n_1,...,n_h|0)
  =(-1)^{d-h}C(g;h|d;n_1,...,n_h|1),
\]
we have
\[
C(g;h|d;n_1,n_2,...,n_h)=0\;\;\; \textup{ if } h>1.
\]
For $g=0$, we have
\[
C(0;1|d;d)=\frac{1}{d^2}.
\]
For $g>0$, we get
\begin{eqnarray*}
C(g;1|d;d)&=&
\int_{{\overM_{g,1}}_{U(1)}}\frac{(-1)^g\lambda_g}{\lambda(\lambda-d\psi_1)}
c_g(\hodge^\vee(\lambda))c_g(\hodge^\vee(-\lambda)),
\end{eqnarray*}
where $\lambda_g=c_g(\hodge)$.  But $c(\hodge)c(\hodge^\vee)=1$
Mumford~\cite[5.4]{M}. It is straightforward to check that this
implies
\[
c_g(\hodge^\vee(\lambda))c_g(\hodge^\vee(-\lambda))=(-1)^g\lambda^{2g}.
\]
It follows that
$$ 
C(g;1|d;d)=d^{2g-2}\int_{\overM_{g,1}}\psi_1^{2g-2}\lambda_g
      =d^{2g-2}b_g. \eqno{\qed}
$$

\bibliographystyle{gtart}
\bibliography{link}

\end{document}